\documentclass[10pt]{article}
\usepackage{latexsym}
\usepackage[all]{xy}
\usepackage{amsmath}
\usepackage{amssymb}
\usepackage{amsfonts}

\newtheorem{thm}{Theorem}[section]

\newtheorem{df}[thm]{Definition}

\newtheorem{Q}[thm]{Question}

\begin{document}

\title{\textbf{Higher and derived stacks: a global overview}}
\bigskip
\bigskip

\author{\bigskip\\
Bertrand To\"en \\
\small{Laboratoire Emile Picard}\\
\small{UMR CNRS 5580} \\
\small{Universit\'{e} Paul Sabatier, Toulouse}\\
\small{France}}

\date{December 2005}

\maketitle

\begin{abstract}
These are expended notes of my talk at the summer institute in algebraic geometry (Seattle, July-August 2005), whose main 
purpose is to present a global overview on the theory of higher and derived stacks. This text is far from being
exhaustive but is intended to cover 
a rather large part of the subject, starting from the motivations and the foundational material, passing through 
some examples and basic notions, and ending with some more recent developments and open questions. 
\end{abstract}

\medskip

\tableofcontents

\bigskip

\section{Introduction}

The notion of algebraic (1-)stack was introduced in the late sixties, and since then it has been highly developed
and has now become a full theory by its own: it is based on solid foundational material (existence of
a nice 2-category of algebraic stacks, notion of sheaves and derived categories \dots), it 
contains many interesting and geometrically meaningful examples (the stack of stable maps, the stack of coherent sheaves
on a variety \dots), many theories have been developed for stacks themselves (intersection theory, l-adic formalism,
vanishing theorems,  motivic cohomology, Riemann-Roch formula, motivic integration \dots) and these
theories have applications
to several other contexts (Gromov-Witten invariants, birational geometry, arithmetic geometry, Hodge theory \dots). 
I think everyone would agree today that the theory of algebraic stack is an important theory. 

Approximately ten years ago C. Simpson introduced in \cite{s3} a notion of \emph{algebraic n-stack}, and more recently
notions of \emph{derived scheme} and of \emph{derived algebraic n-stacks}
have been introduced in \cite{tv3,hagII,lu}. The purpose of this text is to give an overview on the recent works 
on the theories of higher algebraic stacks and of higher derived algebraic stacks, and to show that 
although these theories are not as developed as the theory of algebraic 1-stacks, they 
are based on solid foundational material, contain interesting and geometrically meaningful examples, and
also have interesting developments and applications. \\

This work is organized in three sections. The first part (\S 2) is devoted to the general theory of higher stacks (which is
used all along this work),
in which I tried to explain the motivations and to give some ideas of the foundations of the theory. As there exist several
ways to motivate the theory of higher stacks I had to make a choice and have decided to take the point of view
of moduli theory, but contemplating the theory from another point of view would maybe emphasis different
motivations. I also had to make a choice concerning the foundations of the theory of higher stacks, as
there also exist several possible approaches. I have decided to use the theory of Segal categories (as it 
seems to me the best model for higher categories available today) but also have tried to 
systematically make a bridge with model category theory which provides another approach to 
higher stacks. 

In the second part  of this work (\S 3)  I discuss higher Artin stacks (or algebraic n-stacks). I first present the basic notions
of the theory, such as the definition of higher Artin stacks and some fundamental notions
such as different properties of morphisms (etale, smooth, flat), some sheaf theory (quasi-coherent and l-adic) \dots. 
Then some examples are discussed, trying to present the most significant ones. Finally I present some
developments and applications of higher Artin stacks. Most of them appear already in the literature but some of them 
are only ideas of possible applications and have not been fully investigated. 
The third part of this work (\S 4)
is devoted to derived Artin stacks. It starts by some motivations, and then follows 
the same presentation as the part on higher Artin stacks. \\

To finish this short introduction I would like to thank C. Simpson, G. Vezzosi and J. Lurie for 
numerous conversations about this subject and from which 
I have learned a lot. I am also grateful to J. Kock and to H.H. Tseng for their comments. 

\newpage

\textbf{Conventions:} The expression \emph{stacks} always refer to the notion 
of $\infty$-stacks (in groupoids). To denote the usual notion of stacks in groupoids
as presented in \cite{lm} we will use the expression \emph{1-stacks}. 

All along this text we assume that the reader knows some basic notions of algebraic (1-)stacks (e.g. as
presented in the first chapters of \cite{lm}), and also has some intuitive knowledge of higher category theory
(see e.g. the introduction of \cite{le2}).

\section{Higher stacks}

The main references for this section are \cite{s1,tv,hagI,lu2}
(see also \cite{bre} for a different approach to 2-stacks).

\subsection{Why higher stacks ?}

Moduli theory is about classification of objects and families of objects. Its fundamental concept is that 
of a moduli problem. A moduli problem is a (contravariant) functor $F$, defined on a certain 
category $C$ of geometric objects (e.g. schemes, smooth manifolds, topological spaces \dots), and
whose value $F(S)$ is a structure which is supposed to classify families of objects parametrized by 
the geometric object $S \in C$. In the ancient times the values $F(S)$ of a moduli problem were
taken to be sets, and thus it was implicitly assumed that objects were classified up to 
equality (two points in a sets are or are not equal). However, many moduli problems 
intend to classify objects non only up to equality but also up to isomorphisms, and it was early
recognized that the existence of objects having non trivial automorphisms makes the
set of isomorphism classes of objects badly behaved. Because of this, many interesting
moduli problems could not be representable by conventional geometric objects such as
schemes, smooth manifolds, topological spaces \dots. The theory of 1-stacks (in groupoids) proposes 
a solution to this problem by enhancing the classical notion of moduli functors from set 
valued functors into groupoid valued 
functors. 

One possible starting point of higher stack theory is the observation that 
there exist natural and interesting moduli problems for which objects are classified up to 
a notion of equivalence which is weaker than the notion of isomorphisms. Typical examples are
complexes of abelian groups (or sheaves of abelian groups) up to quasi-isomorphisms, topological spaces up to 
weak homotopy equivalences, or abelian categories up to equivalences of categories. These moduli problems
naturally arise as functors $F : C^{op} \longrightarrow Cat$, for which $F(S)$ must be thought as the
category of objects parametrized by $S$ and equivalences between them. In these new situations, the values
of the moduli functor $F$ are not sets or groupoids anymore but categories. Moreover, the morphisms in these categories
must be "inverted", or "localized", in some sense in order to truly classify objects up to 
equivalences. There exist well known constructions to "invert" a set of morphisms in a category, 
characterized by  universal problems in a 2-categorical context. For instance, 
the Gabriel-Zisman localization is a solution of a universal problem in the 2-category of
categories. There also exists a localization as a solution to a universal problem in the 2-category of "d\'erivateurs" 
(see \cite{der}).  
However, in the same way that the construction sending a groupoid to its set of objects is badly behaved, 
none of these 2-categorical constructions are well behaved. 
It turns out that the meaningful way to "invert" a set of morphisms
in a category is by stating a universal problem in the context of $\infty$-categories (the precise
meaning of this, which requires to fix a theory of $\infty$-categories, 
will be discussed in the next paragraph, and the motivated reader can also consult \cite{to1} for a general discussion
of the localization problem). In particular, "inverting" the equivalences in our moduli functor $F$ provides
an $\infty$-groupoid valued functor. As a conclusion of this short discussion 
it seems to me important to emphasis the following principle: \\

\textbf{Principle 1:} \emph{As 1-stacks appear as soon as objects must be classified up to isomorphism, 
higher stacks appear as soon as objects must be classified up to a notion of equivalence which is
weaker than the notion of isomorphism. } \\

\subsection{Segal categories as models for higher categories}

From very far away an n-category\footnote{The expression \emph{higher category} will always refer to the
weak notion, we will never consider strict higher categories which are somehow useless for the
purpose of stack theory.} is a structure consisting of a set of objects and 
sets of $i$-morphisms for any $0<i\leq n$, together with various kinds of composition (here and later
the integer n can be infinite). A useful inductive point of view
consists of seeing an n-category $A$ as being some sort of category enriched over (n-1)-categories, i.e. of 
a set of objects $Ob(A)$ and for any two objects $a$ and $b$ an (n-1)-category of morphisms $A(a,b)$, together
with composition $A(a,b)\times A(b,c) \longrightarrow A(a,c)$ which is a morphism of (n-1)-categories
and is associative and unital in some sense (to make precise in which sense the associativity and unity axioms hold
is one of the main problem of higher category theory).  There exists many precise definitions of higher categories, 
and I refer the interested reader to \cite{le} for a list of definitions and references. 
Among $\infty$-categories we will mainly be interested in $(1,\infty)$-categories, which by definition 
are the $\infty$-categories whose $i$-morphisms are invertible for any $i>1$. These sorts of higher categories
are extremely important in many contexts as any $\infty$-category obtained by localization from
a 1-category is automatically a $(1,\infty)$-category.  Another way to say 
that an $\infty$-category $A$ is a $(1,\infty)$-category is by stating that for any two objects
$a$ and $b$ the $\infty$-category of morphisms $A(a,b)$ is an $\infty$-groupoid (i.e. all 
its $i$-morphisms are invertible for any $i>0$). Thus, roughly speaking a $(1,\infty)$-category is
a category enriched over $\infty$-groupoids. Moreover, as 
the theory of $\infty$-groupoids is supposed to be (and is for several definitions of \cite{le}) equivalent to the theory of simplicial sets (through some
infinite fundamental groupoid construction), a $(1,\infty)$-category is more or less the 
same thing as a category enriched over simplicial sets (also called \emph{S-categories}). 
This philosophy will explain our choice
of using $S$-categories and more generally Segal categories as models for $(1,\infty)$-categories. 

An $S$-category $A$ is a category enriched over the category of simplicial sets, i.e. 
consists of the data of a set of objects $Ob(A)$, and for any two objects
$a$ and $b$ in $Ob(A)$ a simplicial set of morphisms $A(a,b)$, and 
composition morphisms $A(a,b)\times A(b,c) \longrightarrow A(a,c)$
which are associative and unital (on the nose). A \emph{Segal category} 
$A$ is a weak form of an $S$-category. It consists of 
a set of objects $Ob(A)$, for any two objects $a$ and $b$ 
a simplicial set of morphisms $A(a,b)$, for any three objects $a$, $b$ and $c$ a diagram
of simplicial sets
$$\xymatrix{
A(a,b,c) \ar[r] \ar[d] & A(a,c) \\
A(a,b)\times A(b,c), & }$$
for which the vertical morphism is a weak equivalence of simplicial sets, as well 
as higher structures that I will not make precise. For precise definitions I refer
to \cite{hs,pe,ber}. The main difference between an $S$-category and a Segal category is
that the composition in a Segal category is only defined up to a weak equivalence.
The $S$-categories are precisely the Segal categories for which the vertical 
morphism above is an isomorphism (as well as similar conditions for the higher structures), 
and thus Segal categories generalize $S$-categories. In fact, the two 
notions are  equivalent in some sense (see \cite{ber}) and the reader can
think only in terms of $S$-categories, keeping in mind that Segal categories behave better
for certain purposes and that using $S$-categories could be rather technical at some point. \\

The theory of Segal categories works in a very similar
manner to usual category theory and most (if not all) of
the standard categorical notions can be reasonably defined in the Segal setting.
Here follows a sample of examples (once again we refer to
the overview \cite{tv} for more details).

\begin{enumerate}

\item \textbf{Categories, S-categories and Segal categories:} Segal categories form a category $SeCat$ for the obvious notion of morphisms.
We will use interchangeably the expression \emph{morphisms between Segal categories} and 
\emph{functors between Segal categories}. 

There is a fully faithful functor $S-Cat \longrightarrow SeCat$ from the
category of $S$-categories into the category of Segal categories, and thus
any $S$-category will be considered as a Segal category. Moreover, 
as there is a fully faithful functor $Cat \longrightarrow S-Cat$, 
we also get a full embedding of the category of categories to the
category of Segal categories, and will consider
categories as a special kind of Segal categories (they are the ones for which 
the simplicial sets of morphisms are discrete). 

\item \textbf{Homotopy categories:} Any Segal category $A$ possesses a homotopy
category $Ho(A)$ (which is a category in the usual sense),
having the same objects as $A$, and for two objects $a$ and $b$ morphisms
in $Ho(A)$ are given by $Ho(A)(a,b)=\pi_{0}(A(a,b))$. The composition is 
induced by the diagram
$$\xymatrix{
\pi_{0}(A(a,b,c)) \ar[r] \ar[d]^-{\simeq} & \pi_{0}(A(a,c)) \\
\pi_{0}(A(a,b))\times \pi_{0}(A(b,c)). & }$$
The fact that this composition is associative follows from the higher structures on $A$. 
The functor $SeCat \longrightarrow Cat$ sending $A$ to $Ho(A)$ is left adjoint 
to the embedding $Cat \longrightarrow SeCat$. 

A morphism (between two objects $a$ and $b$) in a Segal category $A$ is a zero simplex of the simplicial set $A(a,b)$.
A morphism is an \emph{equivalence} in $A$ if its image in $Ho(A)$ is an isomorphism.

\item \textbf{Equivalences of Segal categories:}
For a morphism of Segal categories $f : A \longrightarrow B$, we say that
$f$ is \emph{essentially surjective} (resp. \emph{fully faithful}) if
the induced functor $Ho(f) : Ho(A) \longrightarrow Ho(B)$ is essentially surjective (resp.
if for any two objects $a$ and $b$ in $A$ the induced morphism
$f_{a,b} : A(a,b) \longrightarrow B(f(a),f(b))$ is an
equivalence of simplicial sets). We say that $f$ is an \emph{equivalence} if it is
both fully faithful and essentially surjective. When $A$ and $B$ are categories, 
and thus $Ho(A)=A$ and $Ho(B)=B$, this notion of equivalence is the usual 
notion of equivalences of categories. 

\item \textbf{The model category of Segal categories:} The foundational result about Segal categories is the existence
of a model structure whose weak equivalences are the equivalences above (see \cite{hs,pe,ber}). 
To be precise this model structure does not exist on the category of
Segal categories itself but on a slightly larger category of Segal precategories, but
I will simply neglect this fact.  For this model structure, every object
is cofibrant, but not every Segal category is
a fibrant object, and in general fibrant objects are quite difficult to
describe (see however \cite{ber2}). The existence of this model structure is far from being 
formal and has many consequences. First of all it can be used to state that the
theory of $S$-categories and of Segal categories are equivalent in some sense, 
as it is known that their model categories are Quillen equivalent (see \cite{ber}). 
The model category of Segal categories can be shown to be enriched over itself (i.e.
is an internal model category in the sense of 
\cite[\S 11]{hs}, see also \cite{pe,ber}). This implies that
given two Segal categories $A$ and $B$ it is possible to associate a Segal category
of morphisms
$$\mathbb{R}\underline{Hom}(A,B):=\underline{Hom}(A,RB),$$
where $RB$ is a fibrant model for $B$ and $\underline{Hom}$ denotes the
internal $Hom$'s in the category of Segal categories.
From the point of view of $\infty$-categories, $\mathbb{R}\underline{Hom}(A,B)$
is a model for the $\infty$-category of (lax) functors
from $A$ to $B$. In general, the expression \textit{$f : A \longrightarrow B$ is
a morphism of Segal categories} will mean that $f$ is an object
in $\mathbb{R}\underline{Hom}(A,B)$, or equivalently 
a morphism in the homotopy category $Ho(SeCat)$. 
In other words we implicitly allow
ourselves to first take a fibrant replacement of $B$ before considering
morphisms into $B$.

\item \textbf{The 2-Segal category of Segal categories:} Considering fibrant Segal categories 
and their internal Homs as above 
provides a category enriched over $SeCat$, denoted by
$\mathbf{SeCat}$. This is a \emph{2-Segal category} (see \cite{hs}), and is a model for the
$\infty$-category of $(1,\infty)$-categories. I will not really use 
the 2-Segal category $\mathbf{SeCat}$ in the sequel, but it is a good 
idea to keep in mind that it exists. 

\item \textbf{Segal groupoids and delooping:} There is a notion of \emph{Segal groupoid}: by definition 
a Segal category $A$ is a Segal groupoid if its homotopy category
$Ho(A)$ is a groupoid in the usual sense. 

For any Segal category $A$, we can define
its geometric realization $|A|$, which is the diagonal simplicial
set of the underlying bi-simplicial of $A$ (see \cite[\S 2]{hs}, where
$|A|$ is denoted by $\mathcal{R}_{\geq 0}(A)$). The construction
$A \mapsto |A|$ has a right adjoint, sending a simplicial 
set $X$ to its fundamental Segal groupoid $\Pi_{\infty}(X)$
(it is denoted by $\Pi_{1,se}(X)$ in 
\cite[\S 2]{hs}). 
By definition, the set of objects of $\Pi_{\infty}(X)$
is the set of $0$-simplex in $X$, and 
for two points $(x,y) \in X_{0}^{2}$ the simplicial set
of morphisms $\Pi_{\infty}(X)(x,y)$ is the subsimplicial set of
$X^{\Delta^{1}}$ sending the two vertices of $\Delta^{1}$ to
$x$ and $y$. A fundamental theorem states that 
the constructions $A \mapsto |A|$ and
$X \mapsto \Pi_{\infty}(X)$ provide an equivalence
between the homotopy theories of
Segal groupoids and of simplicial sets (see \cite[\S 6.3]{pe}). This last
equivalence is a higher categorical version of 
the well known equivalence between the homotopy theories of 
$1$-truncated homotopy types and of groupoids.

\item \textbf{Localization of Segal categories:} Given a Segal category $A$ and a set of morphisms $S$ in $Ho(A)$, 
there exists a Segal category $L(A,S)$ obtained by  \emph{inverting the arrows in $S$}. This construction
is the Segal analog of the Gabriel-Zisman localization for categories. By definition,
the Segal category $L(A,S)$ comes with a localization morphism $l : A \longrightarrow L(A,S)$
satisfying the following universal property: for any Segal category $B$, the
induced morphism
$$l^* : \mathbb{R}\underline{Hom}(L(A,S),B) \longrightarrow \mathbb{R}\underline{Hom}(A,B)$$
is fully faithful, and its essential image consists of
morphisms $A \longrightarrow B$ sending morphisms of $S$ into equivalences in $B$ (i.e.
isomorphisms in $Ho(B)$). The fact that $L(A,S)$ always exists is not an easy 
result (see e.g. \cite{hs}, or \cite{to2} for a linear analog). 
When applied to the case where $A$ is a category considered as a Segal category, the construction
$L(A,S)$ described above coincides, up to an equivalence, with the simplicial localization
construction of \cite{dk1}. It is important to note that $Ho(L(A,S))$ is naturally 
equivalent to the Gabriel-Zisman localization $S^{-1}A$, but that in general
the natural morphism
$$L(A,S) \longrightarrow S^{-1}A$$
is not an equivalence (examples will be given below). 

\item \textbf{Model categories and Segal categories:} Given a model category $M$, 
we can construct a Segal category $LM:=L(M,W)$ by localizing
$M$ (in the sense as above) along its subcategory of equivalences $W$.
This provides a lot of examples of Segal categories. Using the main result of
\cite{dk3} the Segal categories $LM$ can be explicitly described in terms
of mapping spaces in $M$. In particular, when $M$ is a simplicial model category
$LM$ is equivalent to $Int(M)$, the $S$-category of fibrant-cofibrant objects in $M$. 
For the model category of simplicial sets we will use the notation $Top:=LSSet$, for which 
one model is the $S$-category of Kan simplicial sets. For a simplicial set 
$X$, considered as an object in $Top$, we have a natural equivalence
$Top(*,X)\simeq X$, showing that $Top$ is not equivalent to 
$Ho(Top)$ the homotopy category of spaces. This is the generic situation, and 
for a general model category $M$ the Segal category $LM$ is not equivalent 
to its homotopy category $Ho(M)=Ho(LM)$. 

The construction $M \mapsto LM$ can be made functorial with respects to 
Quillen functors as follows. For $f : M \longrightarrow N$ a right Quillen functor, 
its restriction to fibrant objects $f : M^{f} \longrightarrow N^{f}$ preserves
equivalences, and thus induces a morphism of Segal categories
$LM^{f} \longrightarrow LN^{f}$. The existence of fibrant replacements implies
that the natural inclusion functors $LM^{f} \longrightarrow LM$ and 
$LN^{f} \longrightarrow LN$ are equivalences. We thus obtain a 
morphism of Segal categories $Lf\footnote{Be careful that 
the "L" in "Lf" stands for "localization" and not for "left derived". In fact, as $f$ is right Quillen
the morphism $Lf$ is a model for the right derived functor $\mathbb{R}f$. In order to 
avoid confusion left derived stuff will be denoted using the symbol $\mathbb{L}$.} : LM \longrightarrow LN$, well
defined in the homotopy category of Segal categories, which is often enough
for applications. 

The functor $Lf : LM \longrightarrow LN$ above can be characterized by a universal property
in the Segal category $\mathbb{R}\underline{Hom}(LM,LN)$, showing that it is
uniquely determined by $f$ and equivalences in $M$ and $N$ (in particular it does not
depend on choices of fibrant replacement functors in $M$ and $N$). For this we consider the 
composite functor $\xymatrix{M \ar[r]^-{f} & N \ar[r]^-{l}& LN}$, as an object 
$l\circ f \in \mathbb{R}\underline{Hom}(M,LN)$. In the same way, we have
$Lf\circ l \in \mathbb{R}\underline{Hom}(M,LN)$. By construction there exists 
a natural morphism $l\circ f \rightarrow Lf$ in the Segal category $\mathbb{R}\underline{Hom}(M,LN)$. 
It is possible to show that $l\circ f \rightarrow Lf$ is initial among morphism
from $l\circ f$ to functors $g : M \longrightarrow LN$ sending equivalences in $M$ to equivalences 
in $LN$. In other words, the functor $Lf$ is the \emph{total right derived functor of $f$}, 
in the sense of Segal categories. 

\item \textbf{Classifying spaces of model categories:} For a model category $M$ with subcategory of equivalences $W$, 
we can consider $LM$ as well as its maximal sub-Segal groupoid 
$LM^{int} \subset LM$ defined as the pull-back of Segal categories
$$\xymatrix{
LM^{int} \ar[r] \ar[d] & LM \ar[d] \\
Ho(M)^{int} \ar[r] & Ho(M),}$$
where $Ho(M)^{int}$ is the maximal subgroupoid of $Ho(M)$. As 
$LM^{int}$ is a Segal groupoid it is determined by its geometric realization by the formula
$LM^{int}\simeq \Pi_{\infty}(|LM^{int}|)$. It is possible to show that 
there exists a natural equivalence of simplicial sets 
$|LM^{int}|\simeq |W|$, where $|W|$ is the nerve of the
category $W$ (it is also its geometric realization as a Segal category). Thus, we have
$LM^{int}\simeq \Pi_{\infty}(|W|)$, and the Segal groupoid
$LM^{int}$ is essentially the same thing as the simplicial set 
$|W|$. This fact explains that topologists often refer to  the simplicial set 
$|W|$ as the \emph{classifying space of objects in $M$}: it truly 
is a model for the $\infty$-groupoid obtained from $M$ by inverting 
the morphisms in $W$. This fact will be highly used in the construction and the description 
of higher stacks. 

\item \textbf{The Yoneda embedding:}
Given a Segal category $A$ there is a Yoneda embedding
morphism
$$h : A \longrightarrow \mathbb{R}\underline{Hom}(A^{op},Top),$$
which is known to be fully faithful (this is the Segal version of
the Yoneda lemma). Any morphism $A^{op} \longrightarrow Top$ in the essential image of this
morphism is called \textit{representable}. Dually, there is a notion of
\textit{corepresentable morphism}.

\item \textbf{Adjunctions:} Given a morphism of Segal categories $f : A \longrightarrow B$, we say that
$f$ has a right adjoint if there exists a morphism $g : B \longrightarrow A$
and a point $h \in \mathbb{R}\underline{Hom}(A,A)(Id,gf)$,
such that for any two objects $a \in A$ and $b\in B$ the natural
morphism induced by $h$
$$\xymatrix{A(f(a),b) \ar[r]^-{g_{*}} & A(gf(a),g(b)) \ar[r]^-{h^{*}} &  A(a,g(b))}$$
is an equivalence of simplicial sets.
This definition permits to talk about adjunction between Segal categories.
An important fact is that a Quillen adjunction
between model categories
$$f : M \longrightarrow N \qquad M \longleftarrow N : g$$
gives rise to a natural adjunction of Segal categories
$$Lf : LM \longrightarrow LN \qquad LM \longleftarrow LN : Lg.$$

\item \textbf{Limits:} Given two Segal categories $A$ and $I$, we say that $A$ has limits (resp. colimits) along $I$ if the
constant diagram morphism $A \longrightarrow \mathbb{R}\underline{Hom}(I,A)$
has a right adjoint (resp. left adjoint). This gives a notion of
Segal categories having (small) limits (resp. colimits), or finite limits (resp. colimits).
In particular we can talk about fibered and cofibered squares, final and initial
objects, left and right exactness \dots.

\item \textbf{The strictification theorem:} Let $M$ be a cofibrantly generated model category, and 
$C$ a category with a subcategory $S \subset C$. 
We consider the category $M^{C}$ of functors from $C$ to 
$M$, and $M^{(C,S)}$ the subcategory of functors sending 
morphisms in $S$ to equivalences in $M$. The notion of equivalences in $M$
induces a levelwise notion of equivalences in $M^{(C,S)}$. An important theorem, called
the \emph{strictification theorem}, states that there exists a natural equivalence
of Segal categories 
$$L(M^{(C,S)}) \simeq \mathbb{R}\underline{Hom}(L(C,S),LM).$$
A proof can be found in \cite{hs} (see also \cite{to2} in the context of dg-categories that can easily be
translated to the simplicial setting). This theorem is very important as it 
provides a rather good dictionary between constructions in the context of model categories and
constructions in the context of Segal categories. For instance, for a model
category $M$ the existence of homotopy limits and colimits in $M$ implies
that the Segal category $LM$ possesses limits and colimits in the sense of $(12)$.  

Another important consequence of the strictification theorem and of the Yoneda
lemma $(12)$ states that any Segal category $A$ possesses fully faithful 
embedding $A \longrightarrow LM$ for some model category $M$. This remark 
implies that model categories and Segal categories are essentially the same thing, 
and this relation can be made precise by showing that Segal categories of the form
$LM$ for $M$ a cofibrantly generated model category are exactly the
locally presentable Segal categories (i.e. the cocomplete Segal categories having 
a set of small generators, see \cite{s4}). 

Finally, the strictification theorem also possesses a relative version, 
for presheaves of model categories on $C^{op}$ (the absolute version
above being for the constant presheaf with values $M$), but we will not reproduce 
it here. This generalized strictification theorem
is important for stack theory as it allows to describe certain homotopy limits (see e.g. 
\cite[App. B]{hagII})
of Segal categories in terms of model categories, and is often a key statement
to check that something is a stack.

\end{enumerate}

The previous list of facts shows that 
through the construction $M \mapsto LM$, model category theory is 
somehow an approximation of Segal category theory, and thus of the theory
of $(1,\infty)$-catgeories. The main advantage of passing from model 
categories to Segal categories is the existence of the internal
Hom object $\mathbb{R}\underline{Hom}$, as well as a gain of functoriality. 
However, model categories are 1-categorical structures, and thus
it is reasonable to say that model categories are in some sense 
\emph{strict forms} of $(1,\infty)$-categories.
We finish this paragraph by stipulated this as another important principle: \\

\textbf{Principle 2:} \emph{Model categories are strict forms of $(1,\infty)$-categories, 
and model category theory is a strict form of the theory of $(1,\infty)$-catgeories.} \\

This principle is not only a conceptual one, and it can be verified dramatically in practice. Typically, 
general constructions are done using Segal categories as they are more functorial, but 
explicit computations are usually done using model category techniques. I personally like to think that 
chosing a model category which is a model for a given Segal category (i.e.
\emph{strictifying} the Segal category) is very much like
chosing a system of local coordinates on a manifold: the intrinsic object is the Segal category, but
the model category is useful to have hands on it. For an example of application of the principle 2
to the construction of higher stacks see the end of the next section. \\

To finish this part on Segal categories we introduce the following notations for a Segal category $A$ and two
objects $a$ and $b$
$$Map_{A}(a,b):=A(a,b) \qquad [a,b]_{A}:=\pi_{0}(A(a,b))=Ho(A)(a,b).$$
When $A$ is clear from the context we will simply write $Map(a,b)$ for $Map_{A}(a,b)$ and
$[a,b]$ for $[a,b]_{A}$.

\subsection{Higher stacks}

We are now ready to explain what are higher stacks (in groupoids). For this let me first remind
the following characterization of the category $Sh(C)$ of sheaves of sets on a Grothendieck site $C$.
There exists a functor, which is the  Yoneda embedding followed by the 
associated sheaf functor, $h : C \longrightarrow Sh(C)$. This functor can be characterized by a universal
property in the following way. First of all for two categories $A$ and $B$ with colimits we will denote by 
$\underline{Hom}_{c}(A,B)$ the category of functors commuting with colimits ("c" stands for
"continuous"). Also, recall from \cite{dhi} the notion of an hypercovering in $C$ (and noticed that 
an hypercovering in $C$ is not a simplicial object in $C$ but only in presheaves of sets over $C$). Then
the functor $h : C \longrightarrow Sh(C)$ is characterized up to equivalence by the
following properties. 

\begin{itemize}
\item The category $Sh(C)$ has colimits. 
\item For any category with colimits $B$, the induced functor
$$h^{*} : \underline{Hom}_{c}(Sh(C),B) \longrightarrow \underline{Hom}(C,B)$$
is fully faithful and its image consists of functors $F : C \longrightarrow B$ such that 
for any object $X \in C$ and any hypercovering $U_{*} \rightarrow X$ in $C$
the natural morphism in $B$
$$Colim_{[n]\in \Delta^{op}} F(U_{n}) \longrightarrow F(X)$$
is an isomorphism\footnote{
Here, as well as in the sequel, we make an abuse of language which is commonly 
used in the litterature. The hypercovering
$U_{*}$ is not a simplicial object in $C$ as each $U_{n}$ is only 
a disjoint union of representable presheaves over $C$, and must be understood 
as a formal disjoint union of objects in $C$. For $U=\coprod U_{i}$ such a formal disjoint union 
the notation $F(U)$ stands for $\prod F(U_{i})$.}.
\end{itemize}

Such a characterization also exists for the 2-category of 1-stacks, but in the setting of 2-categories. The definition
of the Segal category of stacks on $C$ is simply the Segal anolog of these two properties. 

\begin{df}\label{d1}
Let $C$ be a Grothendieck site. A \emph{Segal category of stacks on $C$} is a Segal category $A$ together with 
a morphism $h : C \longrightarrow A$ such that the following two properties are satisfied. 
\begin{enumerate}
\item The Segal category $A$ has colimits. 
\item For any Segal category with colimits $B$, the induced morphism
$$h^{*} : \mathbb{R}\underline{Hom}_{c}(A,B) \longrightarrow \mathbb{R}\underline{Hom}(C,B)$$
is fully faithful and its image consists of morphisms $F : C \longrightarrow B$ such that 
for any object $X \in C$ and any hypercovering $U_{*} \rightarrow X$ in $C$
the natural morphism in $B$
$$Colim_{[n]\in \Delta^{op}} F(U_{n}) \longrightarrow F(X)$$
is an equivalence (i.e. an isomorphism in $Ho(B)$). 
\end{enumerate}
When it exists a Segal category of stacks over $C$ will be denoted by $\mathbf{St}(C)$. 
\end{df}

A fundamental result states that for any Grothendieck site $C$ a Segal category of stacks on $C$ exists and is unique up 
to equivalence. 
Once enough of the basic categorical constructions are extended to the
Segal category setting  and proved to behave correctly, this theorem is not difficult to prove and is proved in a similar way 
as the corresponding statement for categories of sheaves. We start by considering the Segal category 
of prestacks $\widehat{C}:=\mathbb{R}\underline{Hom}(C^{op},Top)$, and we 
define $\mathbf{St}(C)$ a localization of $\widehat{C}$ in order to invert all the morphisms of the form
$$colim_{[n]\in \Delta^{op}} U_{n} \longrightarrow X$$
for all hypercovering $U_{*} \longrightarrow X$ in $C$. The fact that this satisfies the correct 
universal property follows from the definition of the localization and from the fact that 
the Yoneda embedding $C \longrightarrow \widehat{C}$ induces for any cocomplete Segal category $B$
an equivalence 
$$\mathbb{R}\underline{Hom}_{c}(\widehat{C},B) \simeq \mathbb{R}\underline{Hom}(C,B),$$
(this last equivalece can be proved from the strictification theorem, point $(13)$ of \S 2.2). 

The Segal category of stacks $\mathbf{St}(C)$ possesses of course a model category counter part, which is
extremely useful in practice. By definition, 
we start by the projective model structure on $SPr(C)$, the category of simplicial presheaves on $C$
(equivalences and fibrations are levelwise). We then define the model category $SPr_{\tau}(C)$ of stacks
over $C$ as being the left Bousfield localization of $SPr(C)$ along the set of morphisms
$U_{*} \longrightarrow X$ for any hypercovering $U_{*} \rightarrow X$ (here
$U_{*}$ is considered as simplicial presheaf and thus as an object in $SPr(C)$, and
is a model for $Hocolim_{[n]\in \Delta^{op}} U_{n}$ computed in $SPr(C)$). Using the strictification theorem (point $(13)$ of \S 2.2) it
is possible to prove that there are natural equivalences of Segal categories
$$L(SPr(C))\simeq \widehat{C} \qquad L(SPr_{\tau}(C))\simeq \mathbf{St}(C).$$
This last equivalences explain that stacks are modeled by simplicial presheaves. This important fact 
has been first stressed by C. Simpson in \cite{s1}, and then
has been used by several author (see e.g.  \cite{hol,ja2}). Another model for $\mathbf{St}(C)$ is the model
category of simplicial sheaves, originally introduced by A. Joyal and revisited by J. Jardine (see \cite{joy,ja1}). This last model shows
that the Segal category $\mathbf{St}(C)$ only depends on the topos $Sh(C)$, and not 
of the choice of the site $C$. 

By universal properties there exists a natural morphism 
$\pi_{0} : \mathbf{St}(C) \longrightarrow Sh(C)$ which can be
thought of as a trunction functor. This morphism has a fully faithful right adjoint
$Sh(C) \longrightarrow \mathbf{St}(C)$ identifying $Sh(C)$ with the full sub-Segal category 
of $\mathbf{St}(C)$ consisting of discrete objects (i.e. objects $x$ for which for any other object $y$
the simplicial set $\mathbf{St}(C)(y,x)$ is equivalent to a set). In the same way, the 2-category of 1-stacks
can be seen as a full sub-Segal category of $\mathbf{St}(C)$ consisting of 
1-truncated objects (see \cite{hol,hagI}). More generally, the full sub-Segal category $St_{\leq n}(C)$ of
$\mathbf{St}(C)$ consisting of n-truncated objects is a model for the (n+1)-category of
n-stacks in groupoids, and the inclusion $St_{\leq n}(C) \longrightarrow 
\mathbf{St}(C)$ possesses a right adjoint $t_{\leq n} : \mathbf{St}(C) \longrightarrow St_{\leq n}(C)$ called
the \emph{n-th truncation functor} (of course $t_{\leq 0}$ coincides with $\pi_{0}$ described
above). 

As for categories of sheaves the localization morphism (which must be considered as
the associated stack functor)
$$a : \widehat{C} \longrightarrow \mathbf{St}(C),$$
has a fully faithful right adjoint $i : \mathbf{St}(C) \longrightarrow  \widehat{C}$. Concrete models
for $\mathbf{St}(C)$ and $\widehat{C}$ can then be described as follows. A model for
$\widehat{C}$ is the $S$-category $SPr(C)^{cf}$ of cofibrant and fibrant objects in $SPr(C)$. A model 
for $\mathbf{St}(C)$ is the full sub-$S$-category of $SPr(C)^{cf}$ consisting of functors $F : C^{op} \longrightarrow SSet$
such that for any hypercovering $U_{*} \longrightarrow X$ in $C$ the natural morphism
$$F(X) \longrightarrow Holim_{[n]\in \Delta}F(U_{n})$$
is an equivalence. Even more concrete models for the homotopy categories $Ho(\widehat{C})$ and
$Ho(\mathbf{St}(C))$ are given by the homotopy category of presheaves of simplicial sets on $C$ and its
full subcategory consisting of functors satifying the descent condition above. With this 
models, the functor
$\pi_{0} : Ho(\mathbf{St}(C)) \longrightarrow Sh(C)$
mentioned above simply sends a simplicial presheaf $F$ to its sheaf of connected component
(i.e. the sheaf associated to the presheaf $X \mapsto \pi_{0}(F(X))$). \\

An important fact is that the morphism
$a$ is left exact (i.e. commutes with finite limits). This has many interesting exactness consequences
on the Segal catgeory $\mathbf{St}(C)$, as for instance the existence of internal Hom objects (i.e. existence
of stacks of morphisms). These exactness properties are formally the same as the one satisfied by 
the Segal category $Top$, and can be summarized as Segal category versions of the standard 
Giraud's axioms for Grothendieck topos. The three fundamental properties are (see \cite{hagI,tv}): 

\begin{enumerate}
\item The Segal category $\mathbf{St}(C)$ has colimits and a set of small generators (this implies that 
it also has limits, which can also be seen directly).
\item Sums in $\mathbf{St}(C)$ are disjoints: for any family of objects $\{x_{i}\}_{i\in I}$ in $\mathbf{St}(C)$ and any 
$i_{1}\neq i_{2}$ in $I$ the following
diagrams
$$\xymatrix{
\emptyset \ar[r] \ar[d] & x_{i_{2}} \ar[d] & & & x_{i_{1}} \ar[r] \ar[d] & x_{i_{1}} \ar[d] \\
x_{i_{1}} \ar[r] & \coprod_{i\in I}x_{i} & & & x_{i_{1}} \ar[r] & \coprod_{i\in I}x_{i}}$$
are cartesian.
\item Equivalence relations are effective in $\mathbf{St}(C)$: for any groupoid object
$X_{1} \rightrightarrows X_{0}$ with quotient $|X_{*}|$, the natural morphism
$X_{1} \longrightarrow X_{0}\times_{|X_{*}|}X_{0}$ is an equivalence in $\mathbf{St}(C)$. 
\end{enumerate}

These three properties can be taken as the definition of a \emph{Segal topos}. I refer
to \cite{tv,lu2} for more on this notion. It can be proved that 
Segal topos  are precisely the Segal categories which are exact localizations of 
Segal categories of the form $\mathbb{R}\underline{Hom}(T,Top)$ for some
Segal category $T$. An important remark however is that there exists Segal topos which are
not exact localizations of Segal categories of the form $\mathbb{R}\underline{Hom}(C,Top)$ for some
category $C$, showing that there exists exotic Segal topos (i.e. which are not determined by a
topos in the usual  sense). Such an exotic Segal topos will be used to develop 
the theory of derived stacks later in this paper. Another example is the Segal category
$\mathbf{St}(k)/F$, for a stack $F$ which is not a sheaf, which is a Segal topos not 
generated by a Grothendieck site in general. \\

To finish this section on higher stacks I would like to give one particular example
of principle $2$ of \S 2.2 in action, concerned with the construction of higher stacks from 
model category data. 

Let $C$ be a Grothendieck site. We are looking for a general procedure to construct 
stacks over $C$, i.e. simplicial presheaves with the descent conditions. From the point 
of view of moduli theory, a stack $F$, which is modeled by a simplicial presheaf $F : C^{op} \longrightarrow SSet$, 
represents a moduli problem: for an object $X\in C$, the simplicial set $F(X)$ 
is a classifying space of families of objects over $X$. From the dictionnay between Segal categories and model
categories (see points $(8)$ and $(9)$ of \S 2.2), we can expect $F(X)$ to be the nerve of the subcategory of
equivalences in a model category $M(X)$, depending on $X$, and being a model for the homotopy theory of families of 
objects parametrized by $X$\footnote{In general 
$F(X)$ is only expected  to be a full subsimplicial set (i.e. union of connected components) of the 
nerve of equivalences in $M(X)$ consisting of objects satisfying certain additional conditions (typically
finiteness conditions)}. The starting point is thus a \emph{presheaf of model categories $M$ on $C$}, 
also called a \emph{Quillen presheaf}: it consists for any $X\in C$ of a model category $M(X)$, and for
any morphism $f : X \longrightarrow Y$ of a left Quillen functor $f^{*} : M(Y) \longrightarrow M(X)$
satisfying $f^{*}\circ g^{*}=(g\circ f)^{*}$ (there is of course a dual notion with right Quillen functors). 
From such a Quillen presheaf we construct a prestack sending $X$ to 
$F(X):=|WM(X)^{c}|$, the nerve of equivalences in $M(X)^{c}$ (i.e. between cofibrant objects), and $f : X \longrightarrow Y$ to the
induced morphism $f^{*} : F(Y) \longrightarrow F(X)$. Note that the restriction to cofibrant objects
is necessary to insure that $f^{*}$ preserves equivalences. However, as for any model category 
$N$ the nerve of equivalences in $N$ and in $N^{c}$ are naturally equivalent to each others, 
$F(X)$ is a classifying space of objects in $M(X)$, as required. 

In this way, for any presheaf of model categories $M$ we obtain a simplicial presheaf 
$F$, which by point $(9)$ of \S 2.2 can be see as the $\infty$-prestack 
of objects in $M$ up to equivalences. The next step is to add conditions on $M$ to insure that 
the prestack $F$ is a stack, i.e. satisfies the descent condition for hypercoverings. 
For any hypercovering $U_{*} \longrightarrow X$ in $C$, we can consider the
the cosimplicial diagram of model categories $n \mapsto M(U_{n})$ (we make here the same
abuse of notation as before, as $U_{n}$ is not an object in $C$ but only a formal 
disjoint union of such, and $M(U_{n})$ means the product of the values of $M$ over the
various components of $U_{n}$). We consider $Sect(U_{*},M)$, the category of 
global section of this cosimplicial category: its objects are families of objects $x_{n}\in M(U_{n})$ for any $n$, together with
morphisms $u^{*}(x_{m}) \rightarrow x_{n}$ in $M(U_{n})$ for any 
simplicial map $u : U_{n} \longrightarrow U_{m}$,  satisfying the usual cocycle condition (see \cite[App. B]{hagII}).
There exists a natural Quillen model structure on $Sect(U_{*},M)$ for which the equivalences and
fibrations are defined levelwise. It is then possible to construct a natural adjunction
$$\phi : Ho(M(X)) \longrightarrow Ho(Sect(U_{*},M)) \qquad Ho(M(X)) \longleftarrow Ho(Sect(U_{*},M)) : \psi.$$
We say that $M$ \emph{satisfies homotopical descent} (the reader will notice the 
analogy with usual cohomological descent for complexes of sheaves) if the above adjunction satisfies the following 
two conditions:
\begin{itemize}
\item The functor $\phi : Ho(M(X)) \longrightarrow Ho(Sect(U_{*},M))$ is fully faifthul.
\item An object $x_{*}Ê\in Ho(Sect(U_{*},M))$ is in the essential image of $\phi$
if and only if for any $u : U_{n} \longrightarrow U_{m}$ the induced morphism
$$\mathbb{L}u^{*}(x_{m}) \longrightarrow x_{n}$$
is an isomorphisms in $Ho(M(U_{n}))$.
\end{itemize}

An important consequence of the strictification theorem (see point $(13)$ of \S 2.2 as well 
as \cite[App. B]{hagII}) states that with the notations above, the prestack $F$ is a stack 
if $M$ satisfies homotopical descent. As far as I know this is the most powerful way to
construct examples of stacks, and many of the examples of stacks presented in the sequel
are based on this construction.

\section{Higher Artin stacks}

Let $k$ be a commutative ring and $k-Aff$ the category of affine $k$-schemes endowed
with the faithfuly flat  and quasi-compact topology. The Segal category of 
stacks $\mathbf{St}(k-Aff)$ will simply be denoted by $\mathbf{St}(k)$, and its objects 
called \emph{k-stacks}. The ffqc topology being subcanonical the natural morphism
$k-Aff \longrightarrow \mathbf{St}(k)$ is fully faithful, and we will simply identify
$k-Aff$ with its essential image in $\mathbf{St}(k)$ (so any stack equivalent to an affine
scheme will be called an affine scheme). 

Recall that a model for $\mathbf{St}(k)$ is the model category of
presheaves of simplicial sets with the local model structure as in \cite{ja1,dhi}, and thus that 
objects in $\mathbf{St}(k)$ might be described concretely as functors 
$$F : k-Aff^{op}=k-CAlg \longrightarrow SSet,$$
from the opposite of the category of affine k-schemes or equivalently the
category of commutative k-algebras, such that for any ffqc
hypercovering of affine schemes $U_{*} \longrightarrow X$ the natural morphism
$$F(X) \longrightarrow Holim_{[n]\in \Delta}F(U_{n})$$
is an equivalence (once again we make the abuse of notation, as  $U_{n}$ is only
a formal disjoint union of affine schemes). We will often use this description in terms of 
simplicial presheaves in order to construct explicit objects in the Segal 
category $\mathbf{St}(k)$.

Recall also that we have introduced the following notations
$$Map_{\mathbf{St}(k)}(F,G):=\mathbf{St}(k)(F,G) \qquad [F,G]_{\mathbf{St}(k)}:=\pi_{0}(Map(F,G)).$$
Moreover, the subscrite $\mathbf{St}(k)$ will not be mentioned when there are
no ambiguities. From a model category theory point of view, 
if $F$ and $G$ are represented by simplicial presheaves we have
$$Map(F,G)\simeq \mathbb{R}\underline{Hom}(F,G)=\underline{Hom}(QF,RG),$$
where $\underline{Hom}$ are the natural simplicial Hom's of the category of simplicial presheaves, and
$Q$ and $R$ are cofibrant and fibrant replacement functors inside the category of simplicial 
presheaves endowed with its local projective model structure (see \cite{ja1,dhi})

The main references for higher Artin stacks are \cite{s3} and \cite[\S 2.1]{hagII}. 
The approach of \cite{hagII} uses model categories, and concerning notation
the homotopy category of the Segal category $\mathbf{St}(k)$ is 
denoted by $St(k)=Ho(\mathbf{St}(k))$ in \cite{hagII}. In the sequel we will work with the Segal category
$\mathbf{St}(k)$, but the constructions and statements given below can also
be translated into a model category language and considered in $St(k)$ (e.g.
the fiber product in the Segal category $\mathbf{St}(k)$ corresponds
to the
homotopy fiber products in $St(k)$, denoted in \cite{hagII} by $\times^{h}$). 

\subsection{Basic notions}

Higher Artin stacks will form a certain sub-Segal category
of $\mathbf{St}(k)$ of objects obtained as nice quotients from affine schemes.  The definition of 
an $n$-Artin stack goes by induction on $n$ as follows. 

\begin{itemize}
\item A (-1)-Artin stack is an affine scheme. A morphism
$f : F \longrightarrow G$ between stacks is
\emph{(-1)-representable}, or \emph{affine}, if for any affine scheme $X$ and any morphism
$X \longrightarrow G$ the pull back $F\times_{G}X$ is 
a (-1)-Artin stack. 

\item Let us assume that the notion of (n-1)-Artin stacks has been defined, as well as
the notion of (n-1)-representable morphisms (one also says
\emph{(n-1)-geometric morphisms}) and smooth
(n-1)-representable morphisms. 
\begin{itemize}
\item 
A stack $F$ is an \emph{n-Artin stack}
if there exists a disjoint union of affine schemes $X$ and 
a smooth (n-1)-representable and surjective morphism $X \longrightarrow F$
(here surjective must be understood in a sheaf-like sense, 
that for any affine scheme $Y$, any morphism $Y \longrightarrow F$ factors through $X$
locally on the ffqc on $Y$). Such a morphism $X \longrightarrow F$ is called a
\emph{smooth n-atlas for $F$}.
\item A morphism $f : F \longrightarrow G$ between stacks is called \emph{n-representable}
(or \emph{n-geometric})
if for any affine scheme $X$ and any morphism $X \longrightarrow G$ the pull back
$F \times_{G}X$ is an n-Artin stack.
\item An n-representable morphism $f : F \longrightarrow G$ between stacks is called \emph{smooth}
if for any affine scheme $X$ and any morphism $X \longrightarrow G$, there exists 
a smooth n-atlas $U \longrightarrow F \times_{G}X$ such that the composition
$U \longrightarrow X$ is a smooth morphism of schemes. 
\end{itemize}
\item A stack $F$  which is an n-Artin stack for some n is simply called an \emph{Artin stack}. 
If furthermore $F$ is n-truncated (i.e. its values as a simplicial presheaf are n-truncated
simplicial sets, $\pi_{i}(F(X))=0$ for all $i>n$ and all $X\in k-Aff$) then $F$ is called an 
\emph{Artin n-stack}. In the same way, 
a morphism $f : F \longrightarrow G$ between stacks is called \emph{representable}
(or \emph{geometric})
if it is n-representable for some n. 
\end{itemize}

The reader is warned that there is a small discrepancy for the indices in the notions
of n-Artin stack and Artin n-stack. For example a scheme is always an 
Artin $0$-stack, but is only a $1$-Artin stack. It is a $0$-Artin stack if and only 
if its diagonal is an affine morphism (see \cite{hagII} for more details on this). 
To avoid confusion we will not use the terminology \emph{n-Artin stack}  which has been 
introduced only for the need of the inductive definition, and we will stay with the
notion of \emph{Artin stack} and \emph{Artin n-stack} which are the pertinent ones for our purpose. \\

Most of the very basic properties of Artin 1-stacks can be shown to extend to the case of Artin stacks. 
Here follows a sample of results.

\begin{enumerate}

\item \textbf{Properties of morphisms:} Any property P of morphisms of schemes 
which is local for the smooth topology 
extends naturally to a property P of morphisms between Artin stacks (see \cite[1.3.6]{hagII}). This
provides notions of unramified, smooth, etale and flat morphisms. A morphism 
of Artin stacks $F \longrightarrow G$ is an open (resp. closed) immersion if for any 
affine scheme $X$ and any morphism $X \longrightarrow G$ the stack 
$F\times_{G}X$ is a scheme and the induced morphism $F\times_{G}X \longrightarrow X$ is
an open (resp. closed) immersion. 

An Artin stack $F$ is \emph{quasi-compact} if it can be covered by an affine scheme (i.e.
there exists a surjective morphism of stacks $X  \longrightarrow F$ with
$X$ affine). A morphism $f : F \longrightarrow G$ between Artin stacks
is \emph{quasi-compact} if for any affine scheme $X$ and 
any morphisms $X \longrightarrow G$ the stack $F \times_{G}X$ is quasi-compact. 
Finally, by induction on $n$, we say that an n-geometric stack $F$ is 
\emph{strongly quasi-compact} if it is quasi-compact and if the diagonal
$F\longrightarrow F\times F$ is strongly quasi-compact (with the convention  that 
a $0$-geometric stack is strongly quasi-compact if it is quasi-compact). 

Finally, an Artin stack $F$ is \emph{locally of finite presentation} if it has a
smooth atlas $U \longrightarrow F$ such that the scheme $U$ is locally of
finite presentation. An Artin stack $F$ is \emph{strongly of finite presentation} 
if it is strongly quasi-compact and locally of finite presentation.

\item \textbf{Presentation as quotient stacks:} The full sub-Segal category of $\mathbf{St}(k)$ consisting of Artin stacks is stable by 
finite limits and disjoint unions. Moreover, a stack $F$ is an Artin n-stack if and only it is equivalent 
to the quotient stack of  a groupoid object $X_{1} \rightrightarrows X_{0}$ with 
$X_{0}$ and $X_{1}$ being Artin (n-1)-stacks and the morphisms $X_{1} \longrightarrow X_{0}$ being smooth
(see \cite[\S 1.3.4]{hagII}). As in the usual case of Artin 1-stacks, the geometry of such a  quotient stack
is the equivariant geometry of the groupoid $X_{1} \rightrightarrows X_{0}$. This also provides a systematic way
to construct examples of higher stacks by taking quotient of schemes by Artin group stacks. For instance,
the quotient stack of a scheme by an action of an Artin group 1-stack is in general 
an Artin 2-stack. 

\item \textbf{Gerbes:} An Artin stack $F$ is a \emph{gerbe} if its 0-truncation $\pi_{0}(F)$ 
is an algebraic space and if the natural morphism $F \longrightarrow \pi_{0}(F)$ is
flat. It can be shown that an Artin stack $F$ is a gerbe if and only if 
the projection 
$$I_{F}:=F\times_{F\times F}F \longrightarrow F$$
is flat ($I_{F}$ is called the \emph{inertia stack} of $F$, and is the stack 
of morphisms from the constant stack $S^{1}:=K(\mathbb{Z},1)$ to $F$) (see \cite{to3}). 

By  generic flatness it can thus been shown that any Artin n-stack $F$ strongly of finite presentation
over $Spec\, k$ possesses a finite decreasing sequence of closed
substacks
$$\xymatrix{F_{r}=\emptyset \ar[r] & F_{r-1} \ar[r] & \dots & F_{1} \ar[r] & F_{0}=F,}$$
such that each stack $F_{i}-F_{i+1}$ is a gerbe. 

\item \textbf{Deligne-Mumford stacks:} An Artin stack $F$ is a \emph{Deligne-Mumford stack} if 
there exists a smooth altas $U \longrightarrow F$ which is an etale morphism. This is equivalent to 
the fact that the diagonal morphism $F \longrightarrow F\times F$ is unramified. However, the
notion of a Deligne-Mumford n-stack is not very interesting for $n>1$, as the 1-truncation 
$\tau_{\leq 1}F$ is always a Deligne-Mumford 1-stack and the natural morphism 
$F \longrightarrow \tau_{\leq 1}F$ is an etale morphism. Indeed, we can write 
$F$ as the quotient of a groupoid object $X_{1} \rightrightarrows X_{0}$ where
$X_{0}$ is a scheme, $X_{1}$ is a Deligne-Mumford (n-1)-stack and the morphism
$X_{1} \longrightarrow X_{0}$ is etale. The morphism $X_{1} \longrightarrow X_{0}$ being etale, 
it is easy to check that the 0-truncation $\pi_{0}(X_{1})$ is an algebraic space etale over $X_{0}$, and that 
furthermore $\pi_{0}(X_{1}) \rightrightarrows X_{0}$ defines an etale groupoid whose
quotient is equivalent to $\tau_{\leq 1}F$. This shows that $\tau_{\leq 1}F$ is a Deligne-Mumford 1-stack, and
the diagram $X_{0} \longrightarrow F \longrightarrow \tau_{\leq 1}F$ shows that 
the projection $F \longrightarrow \tau_{\leq 1}F$ is etale. 

\item \textbf{Flat and smooth atlases:} A stack $F$ is an Artin stack if and only there exists a scheme $X$ and a
a faithfully flat and locally finitely presented representable morphism $p : X \longrightarrow F$. This means that 
we would not gain anything by defining a generalized notion of being an 
Artin stack by only requiring the existence of a flat atlas. 
To show this we define $\mathcal{X}$ to be the stack of
quasi-sections of the morphism $p$,  as follows: for an affine scheme $S$
a morphism $S \longrightarrow \mathcal{X}$ is by definition given by a commutative
diagram in $\mathbf{St}(k)$
$$\xymatrix{
S' \ar[d]_-{f} \ar[r] & X \ar[d] \\
S \ar[r] & F,}$$
where $f$ is a finite flat morphism. The stack $\mathcal{X}$ can be seen to 
be an Artin stack together with a natural projection $\mathcal{X} \longrightarrow F$. We denote
by $\mathcal{X}^{lci}$ the open substack of consisting of points 
for which  the relative cotangent complex of the morphism $S' \longrightarrow X$ 
is perfect of amplitude contained in $[-1,0]$. Then, an argument of obstruction theory shows that 
the morphism $\mathcal{X}^{lci} \longrightarrow F$ is smooth (note that it is 
automatically representable because of the condition on the diagonal of $F$). Finally, the
morphism $\mathcal{X}^{lci} \longrightarrow F$ is surjective as it is so on points with values in 
algebraically closed fields (because any scheme locally of finite type over an algebraically closed
field contains a point which is Cohen MacCauley).

\item \textbf{Homotopy groups schemes:} For any affine scheme and any morphism $X \longrightarrow F$ of 
stacks we define
the \emph{loop stack at $x$} to be 
$$\Omega_{x} F := X\times_{F} X,$$
which is a stack over $X$. The natural morphism
$$\Omega_{x} F \times_{X} \Omega_{x} F \simeq X \times_{F} X \times_{F} X \longrightarrow X\times_{F}X$$
makes $ \Omega_{x} F$ into a group stack over $X$. 
The $n$-th iterated loop stack is defined by induction 
$$\Omega_{x}^{(n)} F := \Omega_{x}(\Omega_{x}^{(n-1)} F),$$
which is again a group stack over $X$. 
The n-th homotopy sheaf of $F$ at the point $x$ is defined to be
$$\pi_{n}(F,x):=\pi_{0}(\Omega_{x}^{(n)} F)$$
and is a sheaf of groups (abelian for $n>1$) on $X$. 

It can be shown that if $F$ is an Artin stack strongly of finite presentation then 
for any $k$-field $K$ and any morphism $x : Spec\, K \longrightarrow F$, the sheaf 
of groups $\pi_{n}(F,x)$ is representable by a group scheme of finite presentation over $Spec\, K$
(see \cite{to3}). The group scheme $\pi_{n}(F,x)$ is the group of $n$-automorphisms of the points $x$ in $F$, and
are higher analogs of the isotropy groups of Artin 1-stacks. 

\item \textbf{Derived categories of $\mathcal{O}$-modules:} 
To each stack $F$ we can associate a Segal topos $\mathbf{St}(F):=\mathbf{St}(k)/F$ of stacks over $F$. 
The Segal topos possesses a natural ring object $\mathcal{O}_{F}:=\mathbb{A}^{1}\times F \longrightarrow F$ 
making
it into a \emph{ringed Segal topos} (i.e. a Segal topos $\mathbf{St}(F)$ together with a colimit commuting
morphism from  $\mathbf{St}(F)$ to the opposite category of commutative rings). As to any ringed topos
$(T,\mathcal{O})$
is associated a derived category $D(T,\mathcal{O})$ of (unbounded) complexes  $\mathcal{O}$-modules, 
the ringed Segal topos $(\mathbf{St}(F),\mathcal{O}_{F})$ gives rise to a \emph{derived Segal category} $L(F,\mathcal{O}_{F})$ 
of $\mathcal{O}_{F}$-modules (see \cite{tvv}). Its homotopy category will be denoted by 
$D(F,\mathcal{O}_{F}):=Ho(L(F,\mathcal{O}_{F}))$ and is called the derived category 
of $F$. 

However, the notion of derived category of a ringed Segal topos
is bit beyond the scope of this overview and we will rather give an explicit description 
of $L(F,\mathcal{O}_{F})$ as follows. The stack $F$ can be writen as a colimit (in $\mathbf{St}(k)$) of 
a simplicial scheme $X_{*}$. For any $n$, the category $C(X,\mathcal{O}_{X_{n}})$ of complexes of 
(big) $\mathcal{O}_{X_{n}}$-modules on $X_{n}$ can be endowed with a
cofibrantly model category structure for which the equivalences are the quasi-isomorphisms (see \cite{ho2}). 
Moreover, for any transition morphism $a : X_{n} \longrightarrow X_{m}$ the adjunction
$$a^{*} : C(X,\mathcal{O}_{X_{m}}) \longrightarrow C(X,\mathcal{O}_{X_{n}}) \qquad
C(X,\mathcal{O}_{X_{m}}) \longleftarrow C(X,\mathcal{O}_{X_{n}}) : a_{*}$$
is a Quillen adjunction. Passing to the localizations (in the sense of point $(7)$ of \S 2.2) 
of the subcategories of cofibrant objects we obtain
a cosimplicial diagram of Segal categories 
$$\begin{array}{ccc}
\Delta & \longrightarrow & SeCat \\
 n & \mapsto & L(C(X,\mathcal{O}_{X_{n}})^{c}).
\end{array}$$
The homotopy limit of this diagram, taken in the homotopy theory of Segal categories gives 
a Segal category $L(F,\mathcal{O})$ which is the derived Segal category of $\mathcal{O}_{F}$-modules. 
Its homotopy category is by definition $D(F,\mathcal{O}_{F})$. The Segal category 
$L(F,\mathcal{O}_{F})$ is \emph{stable} in the sense of \cite{tv2,to1}, and thus
its homotopy category $D(F,\mathcal{O}_{F})$ inherits a natural triangulated structure. 

It is important to note that $D(F,\mathcal{O}_{F})$ is in general not the derived category 
of a ringed topos (as soon as $F$ is not $0$-truncated, i.e. a sheaf of sets), and this is directly related to the fact that the topos
$\mathbf{St}/F$ is in general not generated by a Grothendieck topos. The derived category 
$D(F,\mathcal{O}_{F})$ can also be identified with the full sub-category of the 
derived category of the simplicial scheme $X_{*}$ consisting of objects satisfying the cohomological descent
condition.  
When $F$ is a scheme, $L(F,\mathcal{O}_{F})$ is the Segal category of 
complexes of big $\mathcal{O}_{F}$-modules, and thus
$D(F,\mathcal{O}_{F})$ is the usual derived category of 
sheaves of $\mathcal{O}_{F}$ on the big ffqc site of $F$.

We define a full sub-Segal category $L_{qcoh}(F)$ of 
$L(F,\mathcal{O}_{F})$ consisiting of 
objects $E$ such that for any affine scheme $X$ and any morphism
$u : X \longrightarrow F$, the object $u^{*}(E)\in D(X,\mathcal{O}_{X})$ 
is a quasi-coherent complex. Objects in 
$L_{qcoh}(F)$ will be called quasi-coherent complexes of $\mathcal{O}_{F}$-modules. When 
$F$ is an Artin stack, it is possible to define a t-structure on $L_{qcoh}(F)$, by 
defining objects with non positive amplitude to be $E\in L_{qcoh}(F)$ such that 
for any affine scheme $X$ and any morphism $u : X\longrightarrow F$, the
complex $u^{*}(E)$ has no non zero positive cohomology sheaves. Dually, 
an object $E$ is of non negative amplitude if for any affine scheme and flat morphism
$u : X \longrightarrow F$, the quasi-coherent complex $u^{*}(E)$ on $X$ has
no non zero negative cohomology sheaves (as sheaves on the small Zariski site
of $X$). The heart of this t-structure is denoted by $QCoh(F)$, and is called
the abelian category of quasi-coherent sheaves on $F$. 

For any morphism of stacks $f : F \longrightarrow F'$ there exists an adjunction
of Segal categories
$$f^{*} : L(F',\mathcal{O}_{F'}) \longrightarrow L(F,\mathcal{O}_{F}) \qquad
L(F',\mathcal{O}_{F'}) \longleftarrow L(F,\mathcal{O}_{F}) : f_{*}
$$
The functor $f^{*}$ preserves quasi-coherent complexes, and induces a functor
$$f^{*} : L_{qcoh}(F') \longrightarrow L_{qcoh}(F).$$
It can be shown that this functor admits a right adjoint $f_{*}^{qoch}$
$$f_{*}^{qcoh} : L_{qcoh}(F') \longrightarrow L_{qcoh}(F).$$
However, in general $f_{*}$ does not preserves quasi-coherent complexes and thus $f_{*}^{qcoh}$ is not induced by 
the functor $f_{*}$ in general. However, 
if $f : F \longrightarrow F'$ is a strongly quasi-compact morphisms between Artin stacks, and
if $E\ L_{qcoh}(F)$ is bounded below, then there exists a natural equivalence in $L(F',\mathcal{O}_{F'})$
between 
$f_{*}(E)$ and $f_{*}^{qcoh}(E)$. 

\item \textbf{The l-adic formalism:} Let $l$ be a number invertible
in $k$. 
For any Artin stack 
$F$, we consider $Et/F$ the full sub-Segal category 
of $\mathbf{St}(k)/F$ consisting of morphisms $u : F' \longrightarrow F$
with $F'$ an Artin stack and $u$ an etale morphism. The Segal category 
$Et/F$ possesses a natural topology induced from the one on 
$\mathbf{St}(k)/F$, and is thus a Segal site (see \cite{tv}). The Segal 
category of stacks over $Et/F$ will be denoted by 
$\mathbf{St}(F_{et})$, and is called the \emph{small etale topos of $F$}. 

The constant sheaf of rings $\mathbb{Z}/l^{i}$ on $Et/F$ endows 
$\mathbf{St}(F_{et})$ with a structure of a ringed Segal topos. The 
derived Segal category of this ringed Segal topos will be denoted
by $L(F_{et},\mathbb{Z}/l^{i})$ (see \cite{tvv}). Once again, the notion 
of a derived Segal category of ringed Segal topos is outside
of the scope of this overview, so it is preferable to give the
following more explicit description of $L(F_{et},\mathbb{Z}/l^{i})$. 
We write $F$ as the colimit in $\mathbf{St}(k)$ of a simplicial diagram 
of schemes $X_{*}$. For each $n$, we consider
the category $C((X_{n})_{et},\mathbb{Z}/l^{i})$ of (unbounded) complexes of 
sheaves of $\mathbb{Z}/l^{i}$-modules
on the small etale site of $X$. The localization of
$C((X_{n})_{et},\mathbb{Z}/l^{i})$ along the quasi-isomorphisms is by definition 
the Segal category $L((X_{n})_{et},\mathbb{Z}/l^{i})$. For each simplicial 
morphism $X_{n} \longrightarrow X_{m}$ there is a natural pull back morphism
$$L((X_{m})_{et},\mathbb{Z}/l^{i}) \longrightarrow L((X_{n})_{et},\mathbb{Z}/l^{i}).$$
We get that way a cosimplicial diagram $n \mapsto L((X_{n})_{et},\mathbb{Z}/l^{i})$ of 
Segal categories  and we set 
$$L(F_{et},\mathbb{Z}/l^{i}) =Holim_{n\in \Delta}
L((X_{n})_{et},\mathbb{Z}/l^{i}),$$
where the homotopy limit in taken in the model category of Segal categories. Finally, 
the natural morphisms $\mathbb{Z}/l^{i} \longrightarrow \mathbb{Z}/l^{i-1}$ induce
natural morphisms of Segal categories
$$L(F_{et},\mathbb{Z}/l^{i}) \longrightarrow L(F_{et},\mathbb{Z}/l^{i-1}),$$
and by definition the l-adic derived Segal category of $F$ is
$$L(F_{et},\mathbb{Z}_{l})=Holim_{i}L(F_{et},\mathbb{Z}/l^{i}).$$

As for the case of complexes of $\mathcal{O}_{F}$-modules, the associated homotopy category 
$D(F_{et},\mathbb{Z}/l^{i}):=Ho(L(F_{et},\mathbb{Z}/l^{i}))$ is not the derived category of 
a Grothendieck topos, and this related to the fact that $Et/F$ is not generated by a Grothendieck site except when
$F$ is an algebraic space.

For any morphism of Artin stacks $f : F \longrightarrow F'$ there exists 
a natural adjunction 
$$f^{*} : L(F'_{et},\mathbb{Z}_{l}) \longrightarrow L(F_{et},\mathbb{Z}_{l}) \qquad L(F'_{et},\mathbb{Z}_{l}) \longleftarrow L(F_{et},\mathbb{Z}_{l}) :
f_{*}.$$
It is also possible to define a direct image with compact supports
$$f_{!} : L(F_{et},\mathbb{Z}_{l}) \longrightarrow L(F'_{et},\mathbb{Z}_{l}),$$
at least when the morphism $f$ is strongly of finite type. This morphism 
has a right adjoint
$$f^{!} : L(F'_{et},\mathbb{Z}_{l}) \longrightarrow L(F_{et},\mathbb{Z}_{l}).$$
These four operations can be completed into six operations by 
introducing a tensor product and a corresponding internal
Hom operations. The six operations can then be used to 
prove a base change formula as well as a trace formula for certain 
kind of l-adic complexes satisfying some finiteness conditions. 
These results are out of the scope of the present overview, and the reader
will find the details in the forthcoming work \cite{tvv}. I would also like 
to mention \cite{be} and \cite{laol} where the l-adic formalism has been
studied for Artin 1-stacks, and \cite{to3} for a particular case of the trace formula for
special Artin stacks (see below for the definition). 

\item \textbf{Tangent and cotangent spaces:} Let $F$ be an Artin stack, 
$X=Spec\, A$ an affine scheme and $x : X \longrightarrow F$ be a morphism of stacks. 
We define a morphism of Segal categories
$$\mathbf{Der}_{F}(X,-) : A-Mod \longrightarrow Top$$
in the following way. For an $A$-module $M$ we
consider the trivial square zero extension $A\oplus M$, and the natural
closed embedding of affine schemes $X \longrightarrow X[M]:=Spec\, A\oplus M$. 
We set 
$$\mathbf{Der}_{F}(X,M):=Map_{X/\mathbf{St}(k)}(X[M],F).$$
It can be shown that there exists a unique objects $\Omega_{F,x}^{1}\in D^{\geq 0}(A)$, in the
positive derived category of $A$-modules with natural equivalences
$$\mathbf{Der}_{F}(X,M)\simeq Map(\Omega_{F,x}^{1},M),$$
where the mapping space on the right hand side is taken in the model category of
unbounded complexes of $A$-modules (see \cite{ho}).
The complex $\Omega_{F,x}^{1}$ is called the \emph{cotangent space of $F$ at $x$}, thought its
not an $A$-module but only a complex of $A$-modules. The negative part of the dual complex of $A$-modules is called the
\emph{tangent space of $F$ at $x$} and is denoted by
$$T_{x}F:=\mathbb{R}\underline{Hom}(\Omega_{F,x}^{1},A)_{\leq 0}\in D^{\leq 0}(A).$$
The relation between $\Omega_{F,x}^{1}$ and the tangent stack is the following. We define
a stack $TF \longrightarrow F$ by setting $TF:=\mathbf{Map}(Spec\, k[\epsilon],F)$, where
$k[\epsilon]$ is the k-algebra of dual numbers and where 
$\mathbf{Map}$ are the internal Homs of $\mathbf{St}(k)$ (i.e. the stacks of morphisms). 
For a point $x : X \longrightarrow F$ as above we have a natural equivalence of stacks over $X$
$$X\times_{F}TF \simeq \mathbb{V}(\Omega_{F,x}^{1}),$$
where $\mathbb{V}(\Omega_{F,x}^{1})$ is the linear stack associated to $\Omega_{F,x}^{1}$
as defined in the example $(2)$ of the section \S 3.2 below. 

It is also possible to glue all the complexes $\Omega_{F,x}^{1}$ for $x : X \longrightarrow F$ varying in the
Segal category of smooth morphisms to $F$ and to obtain 
an object $\Omega_{F}\in L_{qcoh}^{\geq 0}(F)$, called the \emph{cotangent sheaf of $F$}, 
thought its not a sheaf but a complexes of sheaves. The negative part of the dual of 
$\Omega_{F}$, as a complex of $\mathcal{O}_{F}$-modules, 
is called \emph{the tangent sheaf of $F$} and is denoted by $T_{F}$. In general
$T_{F}$ is not quasi-coherent anymore (except when $\Omega^{1}_{F}$ is perfect).
There is of course a natural equivalence
of stacks over $F$
$$TF \simeq \mathbb{V}(\Omega_{F}).$$
Finally, in the section on derived stacks (see \S 4) we will see that 
$\Omega_{F}$ is only the truncated version of a \emph{cotangent complex}
encoding important informations about the deformation theory of $F$.

\item \textbf{Complex Artin stacks and analytic stacks:} Assume now that $k=\mathbb{C}$. 
We can define a Segal category $\mathbf{St}(\mathbb{C})^{an}$ of analytic stacks, as 
well as a notion of Artin analytic n-stacks. We start with $Stein$, the site of Stein analytic spaces
endowed with natural transcendent topology. The Segal category of stacks on $Stein$ is denoted
by $\mathbf{St}(\mathbb{C})^{an}$. The notion of Artin n-stacks in $\mathbf{St}(\mathbb{C})^{an}$ is defined
using a straightforward analog of the algebraic notion. 

The analytification functor provides a functor $a : \mathbb{C}-Aff \longrightarrow Stein$, which is
a continuous morphism of sites. It induces an adjunction on the Segal categories of stacks
$$a_{!} : \mathbf{St}(\mathbb{C}) \longrightarrow \mathbf{St}(\mathbb{C})^{an} \qquad
\mathbf{St}(\mathbb{C}) \longleftarrow \mathbf{St}(\mathbb{C})^{an} : a^{*},$$
where on the level of simplicial presheaves the functor $a^{*}$ is defined
by the formula $a^{*}(F)(X):=F(X^{an})$. The functor $a_{!}$ is denoted
by $F \mapsto F^{an}$ and is called the \emph{analytification functor}. Being an inverse image
functor induced from a continuous morphism of sites it commutes with finite limits. 
Moreover, as it sends smooth morphisms between affine $\mathbb{C}$-schemes it is easy 
to check that it preserves Artin n-stacks.

\end{enumerate}

\subsection{Some examples}

\begin{enumerate}

\item \textbf{Eilenberg-MacLane stacks:}
For
a sheaf of abelian groups $A$ (on the site of affine $k$-schemes), we can consider the stack $K(A,n) \in \mathbf{St}(k)$. The stack 
$K(A,n)$ is characterized, up to equivalence in $\mathbf{St}(k)$, by the following universal property: for any affine $k$-scheme 
$X$ there are functorial bijections
$$\pi_{0}(Map(X,K(A,n)))\simeq H^{n}_{ffqc}(X,A).$$
More generally, there exist functorial equivalences of simplicial sets
$$Map(X,K(A,n))\simeq DK(\mathbb{H}_{ffqc}(X,A)),$$
where $\mathbb{H}_{ffqc}(X,A)$ is the complex of cohomology of $X$ with coefficients
in the sheaf $A$, and $DK$ is the Dold-Kan functor from complexes to 
simplicial sets. This implies in particular that we have
$\pi_{i}(Map(X,K(A,n)))\simeq H^{n-i}_{ffqc}(X,A)$.

The Eilenberg-MacLane stacks can be used to define the cohomology groups of
any stack $F\in \mathbf{St}(k)$ with coefficients in the sheaf of abelian groups $A$ by the formula
$$H^{n}(F,A):=\pi_{0}(Map(F,K(A,n))=[F,K(A,n)].$$
This gives a good notion of cohomology for any Artin stacks with coefficients in some
sheaf of abelian groups. Of course, as we use the ffqc topology this is 
ffqc cohomology, and for a scheme $X$ and a sheaf of groups $A$ we have
$H^{n}(X,A)=H^{n}_{ffqc}(X,A)$. 

Finally, when the sheaf of groups $A$ is represented by a an algebraic space which is flat and locally of finite 
presentation over $Spec\, k$, then $K(A,n)$ is an Artin n-stack. In this case the stack 
$K(A,n)$ is moreover smooth over $Spec\, k$, as this can been checked inductively 
on $n$ (the case $n=1$ being treated in \cite{lm}).  \\

\item \textbf{Linear stacks:} Let $F$ be an Artin stack and let $E \in L(F,\mathcal{O}_{F})$ be a 
quasi-coherent complex over $F$. We define a stack $\mathbb{V}(E)$ over $F$ by
$$\begin{array}{cccc}
\mathbb{V}(E) : & \mathbf{St}(k)/F & \longrightarrow & Top \\
 & (f : F' \rightarrow F) & \mapsto & Map_{L(F',\mathcal{O}_{F'})}(f^{*}(E),\mathcal{O}_{F'}).
\end{array}$$
The stack $\mathbb{V}(E)$ is an generalization of the total affine space associated to 
a quasi-coherent sheaf, and is 
called \emph{the linear stack associated to $E$}. By construction, it is characterized by the
following universal property
$$\pi_{0}(Map_{\mathbf{St}(k)/F}(F',\mathbb{V}(E)))\simeq Ext^{0}(f^{*}(E),\mathcal{O}_{F'}),$$
for any $f : F' \longrightarrow F$ in $\mathbf{St}(k)/F$, and where the $Ext^{0}$ is computed in the
derived category of complexes of $\mathcal{O}_{F'}$-modules. 

The stack $\mathbb{V}(E)$ is an Artin stack if 
$E$ is a perfect complex (i.e. its pull-backs to any affine scheme is quasi-isomorphic to
a bounded complex of vector bundles of finite rank), and the morphism
$\mathbb{V}(E) \longrightarrow F$ is then strongly of finite presentation.
If moreover $E$ is perfect with amplitude
contained in $[a,b]$ then $\mathbb{V}(E)$ is an Artin (b+1)-stack. Finally, if $E$ is perfect and of positive amplitude then 
the morphism $\mathbb{V}(E) \longrightarrow F$ is smooth (see \cite{tova}). 

\item \textbf{The stack of abelian categories:} (see \cite{anel}) For any $k' \in k-CAlg$ we consider $k'-Ab$ the category whose objects
are abelian $k'$-linear categories $A$ which are equivalent to 
$B-Mod$, for some associative $k'$-algebra $B$ which is 
projective and of finite type as a module over $k'$. The morphisms in 
$k'-Ab$ are taken to be the $k'$-linear equivalences. For a morphism
$k' \rightarrow k''$ in $k-CAlg$, there exists a base change functor
$$k'-Ab \longrightarrow k''-Ab$$
sending a category $A$ to the category of $k''$-modules in $A$. This
defines a presheaf of categories on $k-Aff$, and passing to the nerve provides a simplicial
presheaf
$$\begin{array}{cccc}
\mathbf{Ab} : & k-CAlg & \longrightarrow & SSet \\
 & k' & \mapsto & \mathbf{Ab}(k'):=N(k'-Ab).
\end{array}$$

The homotopy groups of the simplicial set $\mathbf{Ab}(k')$ can been described
explicitely in the following way. The set $\pi_{0}(\mathbf{Ab}(k'))$ is the set of
equivalences classes of abelian $k'$-linear categories in $k'-Ab$. For a given
$A\in \mathbf{Ab}(k')$, the group $\pi_{1}(\mathbf{Ab}(k'),A)$ is naturally isomorphic 
to the group of isomorphisms classes of autoequivalences of $A$. The group
$\pi_{2}(\mathbf{Ab}(k'),A)$ is the group of invertible elements in the center
of $A$ (i.e. the automorphism group of the identity functor of $A$). Finally, 
for any $i>2$ we have $\pi_{i}(\mathbf{Ab}(k'),A)=0$. 

The object $\mathbf{Ab}$ is considered as a simplicial presheaf over $k-Aff$, and thus
as a stack $\mathbf{Ab}\in \mathbf{St}(k)$. The simplicial presheaf 
itself is not a stack, and thus the natural morphism
$$\mathbf{Ab}(k') \longrightarrow Map_{\mathbf{St}(k)}(Spec\, k',\mathbf{Ab})$$
is not an equivalence in general. This is due to the fact that 
there exist non trivial twisted form of abelian categories for the etale topology.  The object 
$\mathbf{Ab}\in \mathbf{St}(k)$ should therefore be truly considered as the associated
stack to the simplicial presheaf described above. 

It has been proved by M. Anel that $\mathbf{Ab}$ is an Artin 2-stack locally of finite presentation 
over $Spec\, k$ (see \cite{anel}).  Moreover, for an abelian k-linear category $A$, considered
as a global point $A \in \mathbf{Ab}(k)$, then the tangent space of $\mathbf{Ab}$ at 
$A$ is given by 
$$T_{A}\mathbf{Ab}\simeq HH(A)[2]_{\leq 0},$$
where $HH(A)$ is the complex of Hochschild cohomology of $A$.
Therefore, the Artin 2-stack $\mathbf{Ab}$ is a global geometric counter part 
of the formal moduli of abelian categories studied in \cite{lvdb,lvdb2}. 

\item \textbf{The stack of perfect complexes:} For any $k'\in k-CAlg$, we consider
$Parf(k')$ the category of flat perfect complexes of $k'$-modules and quasi-isomorphisms
between them. As we restricted to flat complexes for any morphism $k' \rightarrow k''$
there exists a well defined base change functor
$$-\otimes_{k'}k'' : Parf(k') \longrightarrow Parf(k'').$$
Passing to the nerve we get a simplicial presheaf
$$\begin{array}{cccc}
\mathbf{Parf} : & k-CAlg & \longrightarrow & SSet \\
 & k' & \mapsto & N(Parf(k')),
\end{array}$$
that we consider as an object in $\mathbf{St}(k)$. Using the techniques of left 
Quillen presheaves presented at the end of \S 2.3, it is possible to prove that 
the above simplicial presheaf is already a stack, and therefore
that $\mathbf{Parf}(k')$ is equivalent to $Map(Spec\, k',\mathbf{Parf})$, and 
is a classifying space for perfect complexes of $k'$-modules. 
The set $\pi_{0}(\mathbf{Parf}(k'))$ is in natural bijection with the
set of isomorphisms classes of $D_{parf}(k')$, the perfect derived category of $k'$. 
For a given perfect complex $E\in \mathbf{Parf}(k')$, the group
$\pi_{1}(\mathbf{Parf}(k'),E)$ is naturally isomorphic to 
the automorphism group of the object $E \in D_{parf}(k')$. 
Moreover, the higher 
homotopy group
$\pi_{i}(\mathbf{Parf}(k'),E)$ can be identified with $Ext^{1-i}(E,E)$
for any $i>1$. This provides a rather complete understanding of the
stack $\mathbf{Parf}$. 

The stack $\mathbf{Parf}$ is not truncated as it classifies perfect complexes of abitrary amplitude, and thus
can not be an Artin n-stack for any $n$. However, 
it can be written as a union of substacks $\mathbf{Parf}^{[a,b]}$ of complexes of amplitude
contained in $[a,b]$. It is a theorem that the stacks $\mathbf{Parf}^{[a,b]}$ are Artin n-stacks for $n=(b-a+1)$ and
locally of finite presentation over $Spec\, k$ (see \cite{tova}). Moreover, the natural inclusions
$\mathbf{Parf}^{[a,b]} \hookrightarrow \mathbf{Parf}^{[a',b']}$ are Zariski open immersion, 
and therefore the whole stack $\mathbf{Parf}$ is an increasing union of open Artin substacks. 
Such a stack is called \emph{locally geometric}. The tangent space of 
$\mathbf{Parf}$ taken at a perfect complex $E$ is given by
$$T_{E} \mathbf{Parf}\simeq \mathbb{R}\underline{End}(E,E)[1]_{\leq 0}.$$

The stack $\mathbf{Parf}$ can be generalized in the following way. Let $B$ be an associative and unital 
dg-algebra  over $k$. We assume that $B$ is \emph{saturated}, i.e. that it is perfect 
as complex of $k$-modules, but also as a bi-dg-module over itself. Then, 
for any $k' \in k-CAlg$ we define $\mathbf{Parf}_{B}(k')$ to be the nerve of the
category of quasi-isomorphisms between perfect $B\otimes_{k}^{\mathbb{L}}k'$-dg-modules
(see \cite{tova}).
This defines a stack $\mathbf{Parf}_{B} \in \mathbf{St}(k)$. As above, 
the set $\pi_{0}(\mathbf{Parf}_{B}(k'))$ is in natural bijection with the
set of isomorphisms classes of $D_{parf}(B\otimes_{k}^{\mathbb{L}}k')$, the perfect derived category of $B\otimes_{k}^{\mathbb{L}}k'$. 
For a given  $E\in \mathbf{Parf}_{B}(k')$, the group
$\pi_{1}(\mathbf{Parf}_{B}(k'),E)$ is naturally isomorphic to 
the automorphism group of the object $E \in D_{parf}(B\otimes_{k}^{\mathbb{L}}k')$. 
The higher 
homotopy groups
$\pi_{i}(\mathbf{Parf}(k'),E)$ can be identified with $Ext^{1-i}(E,E)$
for any $i>1$, where the $Ext$-groups are computed in the triangulatd category $D_{parf}(B\otimes_{k}^{\mathbb{L}}k')$.
It is useful to note that the stack $\mathbf{Parf}_{B}$ only depends on the
dg-category $T$ of perfect $B$-dg-modules (i.e. is invariant under derived Morita equivalences). Therefore if $T$ is a dg-category equivalent to the
dg-category of perfect $B$-dg-modules for a saturated dg-algebra $B$ we will simply write
$\mathbf{Parf}_{T}$ instead of $\mathbf{Parf}_{B}$. Using the notations of 
\cite{tova}, $\mathbf{Parf}_{T}$ is the truncation of the derived stack $\mathcal{M}_{T}$. 

The stack $\mathbf{Parf}_{B}$ can be proved to be locally geometric (see \cite{tova}). An important consequence of this
theorem is the existence of a locally geometric stack of perfect complexes on a smooth and proper
scheme $X$ over $k$. Indeed, the derived category $D_{qcoh}(X)$ is known to have a compact generator
$E$ (see \cite{bovdb}). Therefore, if we set $B:=\mathbb{R}\underline{End}(E)$, $B$ is a saturared
dg-algebra such that $D_{parf}(X)\simeq  D_{parf}(B)$, and thus $\mathbf{Parf}_{B}$ can be identified
with $\mathbf{Parf}(X)$ the moduli stack of perfect complexes on $X$. 
An important corollary of the geometricity of $\mathbf{Parf}_{B}$
is thus the geometricity of $\mathbf{Parf}(X)$. 

As a remark, the maximal sub-1-stack $\mathbf{Parf}(X)^{1-rig}\subset \mathbf{Parf}(X)$, consisting
of perfect complexes on $X$ with non negative Ext-groups between themselves, 
is easily seen to be an open substack. The stack $\mathbf{Parf}(X)^{1-rig}$
is therefore an Artin 1-stack. The stack $\mathbf{Parf}(X)^{1-rig}$ has previously been 
shown to be an Artin 1-stack by M. Leiblich in \cite{lie}. 

\item \textbf{Mapping stacks:} The Segal category of stacks $\mathbf{St}(k)$ possesses internal Homs: for
any two objects $F$ and $G$ in $\mathbf{St}(k)$, the morphism 
$$\begin{array}{ccc}
\mathbf{St}(k) & \longrightarrow & Top \\
  H & \mapsto & Map(H\times F,G)
\end{array}$$
is representable by an object $\mathbf{Map}(F,G)\in \mathbf{St}(k)$. 

For a smooth and proper scheme $X$, and an Artin n-stack $F$ locally of finite presentation, 
it can be proved that the stack $\mathbf{Map}(X,F)$ is again an Artin n-stack
locally of finite presentation. The proof of this general fact follows from a
generalization of Artin's representability criterion to higher stacks which can be
found in \cite{lu}. In some cases (i.e. for some particular choices of $X$ and/or $F$), 
it can be proved directly that $\mathbf{Map}(X,F)$ is an Artin n-stack. 
This is for instance the case when $F=\mathbf{Parf}^{[a,b]}$ as we mentioned in the last 
example. Also, when $X$ is finite over $Spec\, k$, the geometricity of 
$\mathbf{Map}(X,F)$ can be proved by an explicit construction of an atlas. 

A much easier situation is for $K$ a finite simplicial set (weakly equivalent to 
a finite simplicial set is enough), considered as
a constant simplicial presheaf over $k-Aff$ and thus as an object in $\mathbf{St}(k)$. For any 
Artin n-stack $F$ the stack $\mathbf{Map}(K,F)$ can be written as a finite limit of
the stack $F$ itself, and thus is again an Artin n-stack. When $K$ represents 
the homotopy type of a compact CW complex $X$, then $\mathbf{Map}(K,F)$ 
should be understood as the stack of non-abelian cohomology of $X$ with coefficients in $F$. 
The fact that $\mathbf{Map}(K,F)$ is an Artin n-stack when $F$ is so is in some sense
a generalization of the fact that the 1-stack of local systems of $X$ is an Artin 1-stack. 

\item \textbf{The stack of saturated dg-categories:} Recall from \cite{tova} and
from the point $4$ above the notion of
a saturated dg-category over the ring $k$. They are the dg-categories quasi-equivalent 
to the dg-category of perfect $B$-dg-module for an associative dg-algebra $B$ which 
is perfect as a complex of $k$-module and as $(B\otimes_{k}^{\mathbb{L}}B^{op})$-dg-module. 

For $k'\in k-CAlg$ we consider
$dgCat_{k'}$ the category of small dg-categories over $k'$. There exists a model category 
structure on $dgCat_{k'}$ whose equivalences are the quasi-equivalences (see \cite{tab}). We
consider $dgCat_{k'}^{cof}$ the subcategory of $dgCat_{k'}$ consisting of
cofibrant objects. We set $\mathbf{dgCat}^{sat}(k')$ to be the nerve of the
category of quasi-equivalences between saturated dg-categories over $k'$. This defines a simplicial presheaf
over $k-Aff$ and thus an object $\mathbf{dgCat}^{sat} \in \mathbf{St}(k)$. 

\begin{Q}\label{Q1}
Is the stack $\mathbf{dgCat}^{sat}$ locally geometric ?
\end{Q}

I believe that the answer to this question is positive. For an integer $n>0$, we define
a substack $\mathbf{dgCat}^{sat,n} \subset \mathbf{dgCat}^{sat}$ of dg-categories
$T$ such that $HH^{i}(T)=0$ for all $i\leq -n$ (here $HH(T)$ denotes the Hochschild cohomology of
the dg-category $T$). It follows from the results of \cite{to2} that  the substack $\mathbf{dgCat}^{sat,n}$
is an $(n+2)$-stack. Moreover, as for a given saturated dg-category $T$ the
Hochschild complex $HH(T)$ is perfect, we clearly have that $\mathbf{dgCat}^{sat}$ is the union
of $\mathbf{dgCat}^{sat,n}$. To answer positively the above question it is then enough 
to show that $\mathbf{dgCat}^{sat,n}$ is an Artin $(n+2)$-stack. This can be 
approached for example by a direct application of the Artin's representability criterion, or even 
better by its extension by J. Lurie to the derived case (see \cite{lu}). 
As expected, the tangent complex should be given by
$$T_{T}\mathbf{dgCat}^{sat}\simeq HH(T)[2]_{\leq 0}.$$

\end{enumerate}

\subsection{Some developments}

\begin{enumerate}

\item \textbf{Some representability statements:} Recall from 
the example $(4)$ of \S 3.2  that for any 
saturated dg-category $T$ of the form $Parf(B)$ for a saturated dg-algebra $B$, 
there exists a locally geometric stack $\mathbf{Parf}_{B}$ classifying 
perfect $B$-dg-modules (or equivalently objects in $T$). As the stack 
$\mathbf{Parf}_{B}$ only depends on $T$ and not on $B$ itself we will denote
it by $\mathbf{Parf}_{T}$. 

As a first consequence of the geometricity of $\mathbf{Parf}_{T}$, if $k$ is a field then
the group $aut(T)$ of self-equivalences of $T$ up to homotopy
($aut(T)$ is really a sheaf of groups) can be seen to be representable
by an algebraic group scheme locally of finite type over $k$. Moreover, it can be
shown that this group only has a countable number of connected components and thus
can be written as an extension
$$\xymatrix{1\ar[r] & aut(T)_{e} \ar[r] & aut(T) \ar[r] & \Gamma \ar[r] & 1}$$
where $\Gamma$ is a countable discrete group and $aut(T)_{e}$ is 
a connected algebraic group of finite type over $k$ (see \cite{tova}). 

Another interesting consequence is the existence of an algebraic space
of simple objects in $T$. For this, we consider the open substack 
$\mathbf{Parf}_{T}^{simp} \subset  \mathbf{Parf}_{T}$ consisting of objects
$E$ in $T$ such that 
$$Ext^{i}(E,E)=0 \; \forall \; i<0 \qquad Ext^{0}(E,E)=k,$$
where the Ext-groups are computed in the triangulated category associated to $T$. 
The substack $\mathbf{Parf}_{T}^{simp}$ is an Artin 1-stack which is a gerbe
over an algebraic space $\pi_{0}(\mathbf{Parf}_{T}^{simp})$ denoted
by $M_{T}^{simp}$. This algebraic space $M_{T}^{simp}$ is a coarse moduli 
space for simple objects in $T$. It can be identified with the quotient stack
$$M_{T}^{simp}\simeq [\mathbf{Parf}_{T}^{simp}/K(\mathbb{G}_{m},1)].$$

We now suppose that $k=\mathbb{C}$.  When $T$ is the dg-category 
of perfect complexes on a smooth and proper variety $X$,
the algebraic space $M_{T}^{simp}$ contains $X$ as a closed and open sub-algebraic space. 
Indeed, an embedding $X\hookrightarrow M_{T}^{simp}$ consists
of sending a point $x\in X$ to the class of the skyscraper sheaf $k(x)$. 
Assume now that $T$ is a dg-model for the Fukaya category of a 
Calabi-Yau variety $X$. It is expected that $T$ is saturated, and thus
the algebraic space $M_{T}^{simp}$ is expected to exist. If a mirror
$X'$ of $X$ exists, then by what we have just seen $X'$ is a 
sup-space of $M_{T}^{simp}$. Therefore, it might be tempting to try 
to construct $X'$ has a well chosen sub-space of $M_{T}^{simp}$. 
In order to be able to say exactly which sub-space $X'$ is 
it is needed to have a reasonable stability condition on $T$, and to try 
to define $X'$ as the sub-space of $M_{T}^{simp}$ classifying 
stable simple objects $E$ in $T$ such that $Ext^{*}(E,E)\simeq Sym(Ext^{1}(E,E)[-1])$. 
This approach suggests that the construction of the mirror only depends
on a good understanding of the Fukaya category of $X$ (i.e. showing that 
it is saturated and  constructing a meaningful stability structure on it).
 
\item \textbf{Motivic invariants:} We will say that an Artin stack $F$ is \emph{special} if
it is strongly of finite presentation and if 
for any field $K$ and any point $x : Spec\, K \longrightarrow F$ the
sheaf $\pi_{i}(F,x)$ is represented by an affine group scheme which is 
unipotent when $i>1$. The class of special Artin stacks already contains several interesting
examples, and they seem to be the reasonable coefficients for non-abelian Hodge cohomology (see \cite{s1}). 

We define an abelian group $K(\mathcal{CH}^{sp}(k))$ by taking the quotient of the
free abelian group over equivalence classes of special Artin stacks by the following three relations:

\begin{enumerate}

\item $$[F\coprod F']=[F]+[F']$$

\item Let $f : F \longrightarrow F'$ be a morphism between special Artin stacks, 
such that for any algebraically closed field $K$ the induced morphism 
$Map(Spec\, K,F)\longrightarrow Map(Spec\, K,F')$ is an equivalence. Then we have
$[F]=[F']$. 

\item Let$F_{0}$ be stack which is either an affine scheme, or
$K(\mathbb{G}_{a},n)$ for some $n>0$. 
Let $f : F \longrightarrow F'$ be a morphism between Artin special stacks such that 
for any morphism $X \longrightarrow F'$ with $X$ an affine scheme, 
there exists a Zariski open covering $U \longrightarrow X$ such that 
$F\times_{F}U$ is equivalent as a stack over $U$ to $F_{0}\times U \longrightarrow U$
(we say that $f$ is \emph{a Zariski locally trivial $F_{0}$-fibration}). 
Then $[F]=[F'\times F_{0}]$.

\end{enumerate}

The group $K(\mathcal{CH}^{sp}(k))$ is made into a ring by setting
$[F].[F']:=[F\times F']$. The ring $K(\mathcal{CH}^{sp}(k))$ is called the Grothendieck
ring of special Artin stacks. It receives a natural morphism from the Grothendieck 
ring of varieties $K(\mathcal{V}(k)) \longrightarrow K(\mathcal{CH}^{sp}(k))$. Here we define
$K(\mathcal{V}(k))$ to be the quotient of the free abelian group over isomorphism classes of
schemes of finite type over $Spec\, k$ by the following two relations:
\begin{enumerate}

\item $$[X\coprod Y]=[X]+[Y]$$

\item Let $f : X \longrightarrow Y$ be a morphism between special Artin stacks, 
such that for any algebraically closed field $K$ the induced morphism 
$X(K)\longrightarrow Y(K)$ is an equivalence. Then we have
$[X]=[Y]$. 
\end{enumerate}

This definition of the Grothendieck ring $K(\mathcal{V}(k))$ only differs from the usual 
one in non-zero characteristic. In general our group $K(\mathcal{V}(k))$ is
the quotient of the usual Grothendieck group obtained by also inverting
the purely inseparable morphisms. 

It can be proved, that if $\mathbb{L}=[\mathbb{A}^{1}]$ then the natural inclusion morphism
$$K(\mathcal{V}(k))[\mathbb{L}^{-1},\{(\mathbb{L}^{i}-1)^{-1}\}_{i>0}]
\longrightarrow K(\mathcal{CH}^{sp}(k))[\mathbb{L}^{-1},\{(\mathbb{L}^{i}-1)^{-1}\}_{i>0}]$$
is an isomorphism (see \cite{to3}). As a consequence we obtain that any additive invariant 
for schemes (i.e. an invariant factorizing through the ring $K(\mathcal{V}(k))$) extends
uniquely as an additive invariant of special Artin stacks). It is possible this way 
to define the motivic Euler characteristic $\chi_{mot}(F)$ of any special Artin stack
as a class in the Grothendieck ring of motives (suitably localized).  Taking the Hodge realization 
we obtain a definition of the Hodge numbers for any special Artin stack. Taking 
the l-adic realization we obtain a version of the trace formula expression the number of
rational point of special Artin stacks over finite field in termes of the trace of the
Frobenius acting on some complex of l-adic cohomology with compact supports.

As an example, for any compact CW complex $X$, represented by a finite simplicial set $K$, 
and for any special Artin stack $F$, the stack $\mathbf{Map}(K,F)$ is again 
a special Artin stack. The Hogde numbers of $\mathbf{Map}(K,F)$ provide 
interesting homotopy invariants of $X$, mesuring in some sense the size of the
space of non-abelian cohomology of $X$ with coefficients in $F$. 

\item \textbf{Hall algebras for dg-categories:} We let $F$ be a stack and we assume that 
it is \emph{locally special} in the sense that it is the union of its open 
special Artin sub-stacks. We define a relative Grothendieck ring
$K(\mathcal{CH}^{sp}(F))$ by taking the free abelian group 
over equivalence classes of morphisms $F' \longrightarrow F$ with $F'$ a special Artin stack, and imposing
the same three relations. The fiber product over $F$ makes $K(\mathcal{CH}^{sp}(F))$ into a commutative
ring (without unit unless $F$ itself special and not only locally special).  

For any morphism $f : F \longrightarrow G$ between locally special 
stacks there exists a natural push-forward
$$f_{!} : K(\mathcal{CH}^{sp}(F)) \longrightarrow K(\mathcal{CH}^{sp}(G)),$$
obtained by sending $F' \longrightarrow F$ to the composite with $f$, which is a morphism 
of abelian groups. When $f$ is 
strongly of finite type, then there also exists a pull-back
$$f^{*} : K(\mathcal{CH}^{sp}(G)) \longrightarrow K(\mathcal{CH}^{sp}(F)),$$
sending $F' \longrightarrow G$ to $F'\times_{G}F \longrightarrow F$, which is 
a morphism of rings. The functorialities $f_{!}$ and $f^{*}$ satisfy the
base change formula when this makes sense. 

Let $T$ be a saturated dg-category and $\mathbf{Parf}_{T}$ the stack of
objects in $T$ as presented in example $(4)$ of \S 3.2. The stack $\mathbf{Parf}_{T}$ is 
locally special, and thus we can consider its
Grothendieck group $K(\mathcal{CH}^{sp}(\mathbf{Parf}_{T}))$. We will use the notation
$$\mathcal{H}_{abs}(T):=
K(\mathcal{CH}^{sp}(\mathbf{Parf}_{T})).$$
We also consider the dg-category $T^{(1)}$ of morphisms in $T$, which is
again a saturated dg-category (see \cite{tova}). There exists a diagram of stacks
$$\xymatrix{
\mathbf{Parf}_{T^{(1)}} \ar[r]^-{c} \ar[d]_-{\pi} & \mathbf{Parf}_{T} \\
\mathbf{Parf}_{T}\times \mathbf{Parf}_{T}. & }$$
The morphism $c$ sends a morphism $x \rightarrow y$ in $T$ to the object $y$. 
The morphism $\pi$ sends a morphism $x \rightarrow y$ to the pair
$(x,y/x)$, where $y/x$ is the cone of the morphism. The morphism $\pi$ can be seen to be 
strongly of finite type  and thus we obtain a natural morphism
$$\xymatrix{
K(\mathcal{CH}^{sp}(\mathbf{Parf}_{T}\times \mathbf{Parf}_{T})) \ar[r]^-{c_{!} \circ \pi^{*}} & 
K(\mathcal{CH}^{sp}(\mathbf{Parf}_{T}))},$$
and therefore a multiplication 
$$\mu : \mathcal{H}_{abs}(T)\otimes \mathcal{H}_{abs}(T) \longrightarrow \mathcal{H}_{abs}(T).$$
It can be checked that this multiplication makes $\mathcal{H}_{abs}(T)$ into 
an associative and unital algebra (by the same argument as in \cite{to4}). The algebra
$\mathcal{H}_{abs}(T)$ is called the \emph{absolute Hall algebra of $T$}. 

The algebra $\mathcal{H}_{abs}(T)$ is  a two-fold generalization of the
usual Hall algebra studied in the context of representation theory (see e.g. \cite{de}). First of all 
it is defined for dg-categories instead of abelian categories, and moreover the base ring
$k$ needs not to be a finite field anymore.  But it is also 
defined by geometric methods and is in some sense a universal 
object mapping to several possible incarnations by means of realization functors. As an 
example, if $k=\mathbb{F}_{q}$ is a finite field, then there exists a morphism of algebras (surjective up to torsion)
$$\mathcal{H}_{abs}(T) \longrightarrow \mathcal{DH}(T),$$
where $\mathcal{DH}(T)$ is the derived Hall algebra defined in \cite{to4}. This morphism
simply sends an object $p : F' \longrightarrow \mathbf{Parf}_{T}$ to the function on 
$\mathbf{Parf}_{T}(\mathbb{F}_{q})$ which counts the number of rational points
in the fiber of $p$. When $T$ is a dg-model for the bounded derived category 
of an abelian category $A$, then $\mathcal{DH}(T)$ contains a copy
of the usual Hall algebra of $A$. This explains how 
$\mathcal{H}_{abs}(T)$ is a geometric counter-part of 
$\mathcal{DH}(T)$, and thus how it generalizes usual Hall algebras. An important 
advantage of $\mathcal{H}_{abs}(T)$ compare to $\mathcal{DH}(T)$ is that it is defined
over $\mathbb{Z}$. It is expected that a suitable generalization of the construction
$T \mapsto \mathcal{H}_{abs}(T)$ to the case of 2-periodic dg-categories
(i.e. dg-categories for which the translation functor $x\mapsto x[2]$
comes equiped with an equivalence with the identity) would give 
a direct construction of quantum enveloppping algebras, generalizing the fact that 
Hall algebras can be used to reconstruct the positive nilpotent part of 
quantum envopping algebras (see \cite{de,to4} for more on the subject).

\item \textbf{The Riemann-Hilbert correspondence:} Let $X$ be a smooth and projective 
complex variety. We associate to it two stacks
$X_{B}$ and $X_{DR}$ in $\mathbf{St}(\mathbb{C})$ as follows. The stack 
$X_{B}$ is the constant stack with values $Sing(X(\mathbb{C}))$, the simplicial sets of
singular simplicies of the topological space of complex points of $X$. The stack 
$X_{DR}$ is defined by its functor of points by $X_{DR}(A):=X(A_{red})$ for 
any $A\in \mathbb{C}-CAlg$. We assume that $F$ is a special Artin stack (as
defined above in point $(2)$) which 
is connected (i.e. the sheaf $\pi_{0}(F)$ is isomorphic to $*$). By example $(5)$ of 
\S 3.2 we know that 
the stack $\mathbf{Map}(X_{B},F)$ is an Artin stack strongly of finite type. 
It can also be proved that the stack $\mathbf{Map}(X_{DR},F)$ is
an Artin stack strongly of finite type (e.g. by using a Postnikov decomposition of $F$). 
A version of the Riemann-Hilbert
correspondence states that there exists a natural equivalence of analytic stacks
(see \cite{s1})
$$\phi : \mathbf{Map}(X_{B},F)^{an}\simeq \mathbf{Map}(X_{DR},F)^{an}.$$
This equivalence is the starting point of a theory of higher non-abelian 
Hodge structures: the stack $\mathbf{Map}(X_{DR},F)$ is considered
as the de Rham non-abelian cohomology of $X$ with coefficients in $F$, 
and the morphism $\phi$ as some kind of integral structure on 
it (at least when $F$ is defined over $\mathbb{Z}$). It is then possible to say what are the Hodge and weight filtrations on 
$\mathbf{Map}(X_{DR},F)$, and to state a definition of a non-abelian 
mixted Hodge structure (see \cite{kps}). \\

\item \textbf{Schematic homotopy types:} Let us assume that $k$ is now a field. 
We let $\mathcal{CH}^{sp}(k) \subset \mathbf{St}(k)$ be the full sub-Segal category 
consisting of special Artin stacks (as defined in point $(2)$). For any connected
simplicial set $K$, we consider the functor between Segal categories
$$\mathcal{CH}^{sp}(k) \longrightarrow Top$$
sending $F$ to $Map(K,F)$. This morphism is not corepresentable by a special Artin stack in general, but
it can be proved to be corepresentable by an object $(K\otimes k)^{sch}\in \mathbf{St}(k)$
which is \emph{local} with respect to the set of objects $\mathcal{CH}^{sp}(k)$. In other words, there exists
a morphism of stacks $u : K \longrightarrow (K\otimes k)^{sch}$ satisfying the
following two conditions:
\begin{enumerate}
\item For any $F\in \mathcal{CH}^{sp}(k)$, the induced morphism
$$u^{*} : Map((K\otimes k)^{sch},F) \longrightarrow Map(K,F)$$
is an equivalence.
\item If $f : G \longrightarrow G'$ is a morphism of stacks such that 
for any $F \in \mathcal{CH}^{sp}(k)$, the induced morphism
$$f^{*} : Map(G',F) \longrightarrow Map(G,F)$$
is an equivalence, then the induced morphism
$$f^{*} :  Map(G',(K\otimes k)^{sch}) \longrightarrow Map(G,(K\otimes k)^{sch})$$
is also an equivalence.
\end{enumerate}

The existence of such a morphism $K \longrightarrow (K\otimes k)^{sch}$ can
easily be deduced from the results of \cite{to5}, and the two above properties
characterizes $(K\otimes k)^{sch}$ uniquely as a stack under $K$. 
The stack $(K\otimes k)^{sch}$ is called the \emph{schematization of $K$ over 
$Spec\, k$}, and is somehow an envelope of $K$ with respect to 
the objects of $\mathcal{CH}^{sp}(k)$. 

The stack $(K\otimes k)^{sch}$ can be proved to satisfy the following 
properties.

\begin{enumerate}
\item We have $\pi_{0}((K\otimes k)^{sch})=*$, and for any 
point $x\in K$, the sheaf $\pi_{i}((K\otimes k)^{sch},x)$
is representable by an affine group scheme which is unipotent for $i>1$. 
Therefore, thought $(K\otimes k)^{sch}$ is not an Artin stack (its diagonal
is not locally of finite type in general), 
it is rather close to be a special Artin stack. In fact it can be shown that 
$(K\otimes k)^{sch}$ is a limit of special Artin stacks, and in some
sense it can be considered as a pro-object in $\mathcal{CH}^{sp}(k)$.
\item The affine group scheme $\pi_{1}((K\otimes k)^{sch},x)$ is
isomorphic to the pro-algebraic completion of the discrete group
$\pi_{1}(K,x)$ over the field $k$. Moreover, for a
finite dimension linear representation $V$ of $\pi_{1}((K\otimes k)^{sch},x)$, corresponding
to a local system $L$ of $k$-vector spaces on $K$, we have
$$H^{*}((K\otimes k)^{sch},V) \simeq H^{*}(K,L).$$
\item When $K$ is simply connected and finite (i.e. each 
$K_{n}$ is finite), then there are isomorphisms
$$\pi_{i}((K\otimes k)^{sch})\simeq \pi_{i}(K)\otimes \mathbb{G}_{a} \qquad if \; char(k)=0$$
$$\pi_{i}((K\otimes k)^{sch})\simeq \pi_{i}(K)\otimes \mathbb{Z}_{p} \qquad if \; char(k)=p>0.$$
This shows that in this case $(K\otimes k)^{sch}$ is a model for the 
rational homotopy type when $k=\mathbb{Q}$ and for the p-adic homotopy type
when $k=\mathbb{F}_{p}$.
\end{enumerate}

In \cite{to5} the construction $K\mapsto (K\otimes k)^{sch}$ has been proposed
as a solution to the \emph{schematization problem} stated in \cite{gr}. In \cite{kpt} the schematization
construction over $\mathbb{C}$ has been used in order to give an alternative
to non-abelian Hodge theory. More precisely, for a smooth and projective
complex manifold $X$, we take $K$ to be the simplicial set of 
singular simplicies of the underlying topological space $X^{top}$ of $X$. The 
schematization of $K$ is simply denoted by 
$(X^{top}\otimes \mathbb{C})^{sch}=(K\otimes \mathbb{C})^{sch}$. The main theorem 
of \cite{kpt} states that there exists an action of the discrete group
$\mathbb{C}^{*}$ on the stack $(X^{top}\otimes \mathbb{C})^{sch}$, called
the \emph{Hodge filtration}. This action can be used to recover all previously known
constructions of the Hodge filtration on cohomology, fundamental group and
rational homotopy groups. It is also possible to prove a purity condition
for this action, that have rather strong consequences on the 
stack $(X^{top}\otimes \mathbb{C})^{sch}$ and thus on the homotopy type
of $X^{top}$. New examples of homotopy types which are not realizable
by smooth projective varieties can be constructed that way. I should also
mention \cite{ol1,ol2} in which a crystalline and a p-adic analog of the constructions above 
have been studied.

\item \textbf{The period map to the moduli of dg-categories:} Let $\mathbf{Var}^{smp}$ be 
the stack of smooth and proper schemes over $Spec\, k$. It is
a 1-stack and thus an object in $\mathbf{St}(k)$. The construction sending a 
smooth and proper scheme $X$ to the dg-category $L_{parf}(X)$ of perfect complexes on $X$
induces a morphism of stacks
$$\phi : \mathbf{Var}^{smp} \longrightarrow \mathbf{dgCat}^{sat}.$$
This morphisms factors through the maximal sub-2-stack $\mathbf{dgCat}^{sat,0}$ 
and thus provides a morphism
$$\phi : \mathbf{Var}^{smp} \longrightarrow \mathbf{dgCat}^{sat,0}.$$
When $k$ is a field of characteristic zero, the tangent of the stack $\mathbf{dgCat}^{sat,0}$ 
at the point $\phi(X)$ can be identified with Hochschild cohomology of $X$ shifted by $2$, and thus we have
$$T_{\phi(X)}\mathbf{dgCat}^{sat,0}\simeq \bigoplus_{p,q}H^{p}(X,\wedge^{q}T_{X})[2-p-q].$$
The map induced on the zero-th cohomology of the tangent spaces by $\phi$ is the natural
embedding
$$H^{1}(X,T_{X}) \longrightarrow H^{0}(X,\wedge^{2}T_{X})\oplus H^{1}(X,T_{X}) \oplus 
H^{2}(X,\mathcal{O}_{X}).$$
This suggests that the morphism $\phi$ is somehow unramified, and thus is a local immersion at least locally on 
$\mathbf{Var}^{smp}$ where $\mathbf{Var}^{smp}$ is an Artin 1-stack. In particular we should get that the
fibers of $\phi$ are discrete (this is not really true because of
stacky phenomenon, but anyway). I think this is a possible geometric approach to a conjecture
(attributed to J. Kawamata) stating that given a given triangulated category $T$ there exists a most finitely many 
smooth and projective varieties having $T$ as perfect derived category.

\end{enumerate}

\section{Derived stacks}

The main references for derived stacks are \cite{tv3,hagII,lu}. 

\subsection{Why derived stacks ?}

We suppose that we are given a moduli 
functor
$$F : k-Aff^{op} \longrightarrow SSet,$$
which is represented by a scheme $X$, or even 
an Artin n-stack also denoted by $X$. The classical problem of 
obstruction theory can be stated as follows: given
any surjective morphism $A \longrightarrow A_{0}$ in $k-CAlg$, 
with kernel $I$ such that $I^{2}=0$, study the fibers of the induced morphism
$$F(A) \longrightarrow F(A_{0}).$$
When $F$ is given  by a concrete moduli problem, there
exists  a complex $L\in D(A_{0})$, which is somehow
"natural" (in the psycological sense of the word), and for any 
point $x\in F(A_{0})$ a class $e\in Ext^{1}(L,I)$, such that 
the fiber at $x$ is non-empty if and only if $e=0$. 

The first observation is that the complex $L$ is by no means unique. In fact, 
there are situations for which there exist different possible 
choices for $L$, all of them being "natural" in some sense. Once again, they are
natural only in the psyclolgical sense of the word and are definitely
not natural in any mathematical sense, unless there will not be any choices. 
For instance, a morphism between moduli functors might not
induce morphisms on the corresponding complexes.  
Moreover, 
forgetting the moduli functor $F$ and only keeping the scheme, or Artin stack, $X$, 
there also exists the cotangent complex of $X$ at the point $x$, 
$\mathbb{L}_{X,x}\in D(A_{0})$, and a natural obstruction class
$e\in Ext^{1}(\mathbb{L}_{X,x},I)$, such that 
the fiber at $x$ is non-empty if and only if $e=0$. A striking remark is that 
in practice, when $F$ is a concrete moduli problem, then 
the two objects $L$ and $\mathbb{L}_{X,x}$ are in general not the same. Even more 
striking is the fact that for a concrete moduli problem $F$ the complex
$\mathbb{L}_{X,x}$ is in general very hard (if not impossible) to compute in terms
of $F$, whereas $L$ has very concrete geometrical description. 

Here is a typical example: let $S$ be a smooth and proper scheme over $k$, 
and $F=\mathbf{Vect}(S)$ be the 1-stack of vector bundles on $S$, which
is an Artin 1-stack. For a point $x : X_{0}:=Spec\, A_{0} \longrightarrow \mathbf{Vect}(S)$, 
corresponding to a vector bundle $E$ on $S\times X_{0}$, the natural candidate
for $L$ is the complex $\mathbb{R}\underline{End}(E,E)^{\vee}[-1] \in D(A_{0})$. 
However, $\mathbb{R}\underline{End}(E,E)^{\vee}[-1]$ being perfect and
not of amplitude contained in $[-1,\infty[$ in general, it can not be
the cotangent complex of any Artin stack locally of finite presentation. 
Moreover, the cotangent complex of the stack $\mathbf{Vect}(S)$
at the point $x$ is not known. This example is not a pathology, and 
reflects the general situation. 

What this example, and many other examples, shows is that 
in general the right complex to consider to understand osbtruction theory 
is $L$, not $\mathbb{L}_{X,x}$, but also that $L$ is not the cotangent complex of any 
Artin stack (locally of finite presentation). One purpose of the notion of \emph{derived
Artin stack} is precisely to provide a new geometric context in which 
the complex $L$ is truly the cotangent complex of some geometric object. In this new
context, $L$ being the cotangent complex of some geometric object will be
natural, now in the mathematical sense of the word, and thus obstruction theory
will become unambiguous. \\

\textbf{Principle 3:} \emph{Derived algebraic geometry is a generalization 
of algebraic geometry for which obstruction theory becomes natural.} \\

Of course the price to pay is that the correct moduli space associated 
to $F$ can not be a scheme or an Artin stack anymore, and 
other kind of geometrical objects are needed, called \emph{derived Artin stacks}. 
In order to guess what these are I would like to come back to our example
of a moduli functor $F$ and to the infinitesimal lifting problem. 

The problem is to understand the obstruction class $e\in Ext^{1}(L,I)$ from a geometric point of view. 
For this, recall that for any $A_{0}$-module $M$, if $A_{0}\oplus M$ denotes the trivial square
zero extension of $A_{0}$ by $M$, then the fibers of $F(A_{0}\oplus M) \longrightarrow F(A_{0})$ 
are isomorphic to $[L,M]=Ext^{0}(L,M)$ (we assume here that $F$ is a set valued functor 
for the sake of simplicity). This suggest that $Ext^{1}(L,I)=[L,I[1]]$ should be the fiber of the
morphism $F(A_{0}\oplus I[1]) \longrightarrow F(A_{0})$. Of course $A_{0}\oplus I[1]$ does not make
sense in rings anymore, but can be defined as a commutative dg-algebra, or better as a simplicial 
commutative algebra with $\pi_{0}(A_{0}\oplus I[1])=A_{0}$, $\pi_{1}(A_{0}\oplus I[1])=I$ and
$\pi_{i}(A_{0}\oplus I[1])=0$ for $i>1$. Therefore we already see that the obstruction space 
$Ext^{1}(L,I)$ will have a functorial description in terms of $F$ as soon as $F$ is extended 
from commutative rings to simplicial commutative rings. Moreover, there exists a 
homotopy pull back diagram of simplicial rings
$$\xymatrix{
A \ar[r] \ar[d] & A_{0} \ar[d] \\
A_{0} \ar[r] & A_{0}\oplus I[1],}$$
suggesting that when $F$ is reasonable there exists a pull back diagram
$$\xymatrix{
F(A) \ar[r] \ar[d] & F(A_{0}) \ar[d] \\
F(A_{0}) \ar[r] & F(A_{0}\oplus I[1]).}$$
The obstruction class $e$ is then expected to be the image of the point $x\in F(A_{0})$ in 
$F(A_{0}\oplus I[1])$, which naturally lives in the fiber of the projection to $F(A_{0})$ and thus
in $Ext^{1}(L,I)$. 

We conclude that the obstruction theory of $F$ can be explained as soon as $F$ is extended
to a functor defined on the category of simplicial rings. If such an extension is given, 
we clearly expect that $Ext^{i}(L,I)$ is the fiber of $F(A_{0}\oplus I[i]) \longrightarrow F(A_{0})$. 
This suggest that once an extension $\widetilde{F}$ to simplicial rings is given then $L$ becomes uniquely determined
by $F$, and should be thought of as the cotangent complex of $\widetilde{F}$. The non uniqueness of
$L$ with respect  to $F$ is then related to the non uniqueness of the extension of $F$ to simplicial rings. \\

The conclusion of this small discussion is: as stacks are functors defined on the category of
commutative rings, derived stacks are functors defined on the category of simplicial rings. I like to draw
the following picture, relating sheaves, 1-stacks, higher stacks and derived stacks all together
$$\xymatrix{
k-CAlg \ar[rrd]^-{1-stacks} \ar[rrdd]_-{stacks} \ar[rr]^-{sheaves} \ar[dd]_-{j} & & Set \\
 & & Groupoids \ar[u]_-{\pi_{0}} \\
 sk-CAlg \ar[rr]_-{derived\; stacks} & & SSet. \ar[u]_-{\Pi_{1}}
 }$$
In this picture, $sk-CAlg$ is the category of simplicial objects in $k-CAlg$, $j$ is the natural inclusion 
functor seeing a $k$-algebra as a constant simplicial object, $\pi_{0}$ is the functor sending
a groupoid to its set of isomorphism classes and $\Pi_{1}$ sends a simplicial set to
its fundamental groupoid. 

An important new feature in the theory of derived stacks is that the category $sk-CAlg$ of
commutative simplicial k-algebras has a natural model category structure, and naturally 
the weak equivalences have to be "inverted" or "localized". Therefore, derived stacks should truly be morphisms 
of Segal categories
$$L(sk-CAlg) \longrightarrow Top,$$
and are not modeled by simplicial presheaves on some Grothendieck site anymore. We will see however
that the Segal category $L(sk-CAlg)$ has a natural extension of the usual ffqc topology, and that derived stacks
can then be viewed as stacks on the \emph{Segal site} $(L(sk-CAlg),ffqc)$. 

\subsection{Basic notions}

We start with the category $sk-CAlg$ of simplicial commutative $k$-algebras. 
There is a natural notion of weak equivalences between objects in $sk-CAlg$, defined
as the morphisms inducing weak equivalences on the underlying simplicial sets (by forgetting the
ring structure). The Segal category $L(sk-CAlg)^{op}$, obtained by localizing the equivalences
in $sk-CAlg^{op}$ is defined to be the Segal category of \emph{derived affine schemes} and is denoted by
$$dk-Aff:=L(sk-CAlg)^{op}.$$
As there exists a simplicial model category structure on $sk-CAlg$, for which the equivalences and fibrations are defined
on the underlying simplicial sets, the Segal category 
$dk-Aff$ can be concretely described as $Int(sk-CAlg)^{op}$ (see point $(8)$ of \S 2.2). 

The Segal category of derived pre-stacks is then defined to be
$$\widehat{dk-Aff}:=\mathbb{R}\underline{Hom}(dk-Aff^{op},Top).$$
Using the dictionnary between model categories and Segal categories (see point $(8)$ of \S 2.2) 
the Segal category $\widehat{dk-Aff}$ can be described by the homotopy theory 
of equivalence preserving functors $sk-CAlg \longrightarrow SSet$. More precisely, 
we can define a model category $M$, by first considering the model category 
of functors $SSet^{sk-CAlg}$ endowed with the levelwise projective model structure, and
then define $M$ as the left Bousfield localization of $SSet^{sk-CAlg}$ along the
equivalences in $sk-CAlg$  (see \cite{hagI} for details). We then have a natural equivalence
between $\widehat{dk-Aff}$ and $Int(M)$. 

The next step is to endow $dk-Aff$ with a topology. It can be shown that the natural
notion of a Grothendieck topology on a Segal category $A$ is nothing else than 
a Grothendieck topology on its homotopy category $Ho(A)$ (see \cite{hagI,tv} for a justification). 
For this we wil need the following important definitions. The fact that these definitions are reasonable extensions of the usual notions is explained in 
\cite[\S 2.2.2]{hagII}.

\begin{df}\label{d2}
A morphism $f : A \longrightarrow B$ in $sk-CAlg$ is 
\emph{flat} (resp. \emph{smooth}, resp. \emph{etale}, resp. a \emph{Zariski open immersion})
if it satisfies the following two conditions:
\begin{enumerate}
\item The induced morphism of affine scheme $Spec\, \pi_{0}(B) \longrightarrow Spec\, \pi_{0}(A)$ is flat
(resp. smooth, resp, etale, resp. a Zariski open immersion).
\item For any $i>0$, the natural morphism
$$\pi_{i}(A)\otimes_{\pi_{0}(A)}\pi_{0}(B) \longrightarrow \pi_{i}(B)$$
is an isomorphism. 
\end{enumerate}
A finite family of morphisms $\{f_{i} :  A \longrightarrow B_{i}\}$ in $sk-CAlg$ is 
a \emph{ffqc covering}, if each $f_{i}$ is flat and if the induced morphism of affine schemes
$$\coprod Spec\, \pi_{0}(B_{i}) \longrightarrow Spec\, \pi_{0}(A)$$ 
is surjective.
\end{df}

We now define a Grothendieck topology on $Ho(dk-Aff)=Ho(sk-CAlg)^{op}$ by 
defining a sieve to be a covering sieve if it contains a ffqc covering in the sense of the definition above. 
It can be checked that this defines a topology on $Ho(dk-Aff)$ and thus by definition a topology on 
the Segal category $dk-Aff$. The ffqc topology on $dk-Aff$ induces a notion 
of hypercoverings in $dk-Aff$, and the Segal category of stacks over $dk-Aff$ can then be
defined in the following way.

\begin{df}\label{d3}
The \emph{Segal category of derived stacks} (also called \emph{$D^{-}$-stacks})
is the full sub-Segal category of $\widehat{dk-Aff}$ consisting of morphisms
$$F : dk-Aff^{op} \longrightarrow Top$$
such that for any ffqc hypercovering $U_{*} \longrightarrow X$
in $dk-Aff$, the induced morphism
$$F(X) \longrightarrow Lim_{n\in \Delta} F(U_{n})$$
is an equivalence in $Top$\footnote{We make here the same abuse of notations
as at the beginning of \S 2.3}. It is denoted by  $\mathbf{dSt}(k)$. 
\end{df}

Like in the case of Segal categories of stacks over a Grothendieck site (see \S 2.3), 
the Segal category $\mathbf{dSt}(k)$ can be characterized by a universal property. 
Also, using the dictionnary between Segal categories and model categories (see point $(8)$ of \S 2.2), 
a concrete model for $\mathbf{dSt}(k)$ is the homotopy theory of functors
$$F : sk-CAlg \longrightarrow SSet$$
satisfying the following three properties:
\begin{enumerate}
\item For any equivalence $A \longrightarrow B$ in $sk-CAlg$ the induced morphism
$F(A) \longrightarrow F(B)$ is an equivalence of simplicial sets.
\item  For any coaugmented co-simplicial object $A \longrightarrow B_{*}$
in $sk-CAlg$, which correspond to a ffqc hypercovering in $dk-Aff$, 
the induced morphism
$$F(A) \longrightarrow Holim_{n\in \Delta} F(B_{n})$$
is an equivalence of simplicial sets. 
\item For any finite family of objects $\{A_{i}\}$ in $sk-CAlg$, 
the natural morphism
$$F(\prod A_{i}) \longrightarrow \prod F(A_{i})$$
is an equivalence of simplicial sets. 
\end{enumerate}

The homotopy theory of these functors can be described by a natural model category, 
called the model category of derived stacks and which is denoted by
$D^{-}k-Aff^{\sim,ffqc}$ in \cite[\S 2.2]{hagII}. To construct derived stacks we will often construct explicit objects
in $D^{-}k-Aff^{\sim,ffqc}$ and then consider them as objects in $\mathbf{dSt}(k)$ through the equivalence
$$\mathbf{dSt}(k) \simeq L(D^{-}k-Aff^{\sim,ffqc}).$$

The natural inclusion morphism $\mathbf{dSt}(k) \hookrightarrow \widehat{dk-Aff}$ has 
an exact left adjoint $a : \widehat{dk-Aff} \longrightarrow \mathbf{dSt}(k)$, called
the \emph{associated derived stacks functor}. The exactness of $a$ implies that 
$\mathbf{dSt}(k)$ does have the same exactness properties as the Segal category $Top$ and
that it is a Segal topos (see \cite{tv}). As a consequence it possesses all small limits and colimits, and has internal Homs. 

Moreover, the ffqc topology can be seen to be subcanonical, and thus the Yoneda embedding
provides a fully faithful functor
$$dk-Aff\simeq L(sk-CAlg)^{op} \hookrightarrow \mathbf{dSt}(k).$$
On the level of simplicial commutative k-algebras this functor will be denoted by
$$\mathbb{R}Spec : L(sk-CAlg)^{op} \hookrightarrow \mathbf{dSt}(k).$$
For any $A\in sk-CAlg$, the derived stack 
$\mathbb{R}Spec\, A$ is explicitly given by 
$$\begin{array}{cccc}
\mathbb{R}Spec\, A : & sk-CAlg & \longrightarrow & SSet \\
 & B & \mapsto & Map_{sk-CAlg}(A,B).
\end{array}$$

The natural embedding $i : k-CAlg \hookrightarrow sk-CAlg$, sending a commutative k-algebra to 
the associated constant simplicial object, induces a morphism on Segal categories of stacks
$$t_{0}:=j^{*} : \mathbf{dSt}(k) \longrightarrow \mathbf{St}(k)$$
called the truncation functor (it is not the same as the truncation $t_{\leq 0}$ defined 
from (underived) stacks to sheaves and discussed in \S 2.3). This functor has a left adjoint
$$i:=j_{!} : \mathbf{St}(k) \longrightarrow \mathbf{dSt}(k)$$
which can be shown to be fully faithful.  In particular, any stack can be seen as 
a derived stack. However, the functor $i$ is not compatible with finite limits, and therefore
certain construction (such as fiber products or internal Homs) will not 
preserve stacks inside the Segal category of derived stacks. Because of this it is important to
keep the notation $i(F)$, when a stack $F$ is considered as a derived stack. 

\begin{df}\label{d4-}
Let $F$ be a stack. A \emph{derived enhancement of $F$} is
a derived stack $\widetilde{F}$ together with an equivalence
$t_{0}\widetilde{F}\simeq F$.
\end{df}

Of course, a given stack $F$ has many different derived enhancement, including 
the trivial one $i(F)$. \\

Using the notion of smooth morphism defined in def. \ref{d3}, the notion of n-geometric stack can be
naturally extended to the notion of \emph{n-geometric derived stack}. As this is a formal 
generalization we will not give the precise definition here (the reader can consult 
\cite{hagII} for more details on the general notion of n-geometric stacks in 
various contexts). The two functors $i$ and $t_{0}$ above are compatible with the
geometricity notions in the sense that $i$ sends n-geometric stacks to 
n-geometric derived stacks, and $t_{0}$ sends n-geometric derived stacks to 
n-geometric stacks. Moreover, a stack $F$ is n-geometric if and only if $i(F)$ is
an n-geometric derived stack. Finally, $i$ and $t_{0}$ are also compatible with the 
notions of flat morphisms, smooth morphisms, etale morphisms and open Zariski morphisms.  

\begin{df}\label{d4}
A derived stack is a \emph{derived Artin n-stack} if 
it is an m-geometric derived stack for some m, and if 
$t_{0}(F)$ is an (Artin) n-stack. A \emph{derived Artin stack} is
an Artin n-stack for some n. 
\end{df}
   
Here is a sample of basic notions and results concerning derived Artin stacks. \\

\begin{enumerate}

\item \textbf{Properties of morphisms and presentations by groupoids:} Both points $(1)$ and $(2)$ 
of the general properties stated in \S 3.1 generalize immediatly to the case of derived Artin stacks. Only two remarks have to be made 
concerning unramified morphisms and closed immersions. First of all, the
notion of formally unramified morphisms in the context of derived Artin stacks is equivalent to the
notion of formally etale morphism (see \cite[Prop. 2.2.2.9]{hagII}). Also, the closed immersions of 
derived Artin stacks are not monomorphisms. In fact, a monomorphism
of derived Artin stacks is automatically formally unramified and thus formally etale, which 
explains why it would not be reasonable that closed immersions be monomorphisms. As a consequence
the notion of closed sub-stack does not make very much sense in the derived setting. 

\item \textbf{Truncation:} For any derived Artin stack $F$, the 
adjonction morphim $it_{0}F \longrightarrow F$ is a closed immersion. 
If fact, if $F$ is locally of the form $\mathbb{R}Spec\, A$, for some $A\in sk-CAlg$, then
$it_{0}$ is locally of the form $i(Spec\, \pi_{0}(A))$. Therefore, 
any derived Artin stack $F$ can be thought of as some kind of 
derived thickening of its truncation $it_{0}(F)$. This derived thickening 
truly behaves as a formal thickening, and for instance the small etale sites
of $F$ and $it_{0}(F)$ coincide (see \cite[Cor. 2.2.2.13]{hagII}). According to definition \ref{d4-}, 
a derived enhancement of a stack can then be thought of as the data of
a formal derived thickening. 

\item \textbf{Derived schemes and Deligne-Mumford stacks:} A derived Artin stack 
$F$ is a \emph{derived scheme} (resp. 
a \emph{a derived Deligne-Mumford stack}) if 
there exists a smooth atlas $U \longrightarrow F$ which is
a Zariski open immersion (resp. etale). It can be shown that if
$F$ is a derived Artin stack then $F$ is a derived scheme (resp. a 
derived Deligne-Mumford stack) if and only if its truncation 
$t_{0}$ is a scheme (resp. a Deligne-Mumford stack) in the non derived
sense. 

\item \textbf{Derived categories of $\mathcal{O}$-modules:}

Like in the underived case any derived stack $F$ has a Segal category 
$L(F,\mathcal{O}_{F})$ of (unbounded) complexes of
$\mathcal{O}_{F}$-modules (see \cite{tvv} for more details). First of all 
the Segal topos $\mathbf{dSt}(k)$ has a natural ring object, 
denoted by $\mathcal{O}$ and represented by $\mathbb{A}^{1}$. The object
$\mathcal{O}$ can also be seen as a colimit preserving functor
$$\mathcal{O} : \mathbf{dSt}(k) \longrightarrow L(sk-CAlg),$$
from $\mathbf{dSt}(k)$ to the Segal category of simplicial commutative k-algebras. The
pair $(\mathbf{dSt}(k),\mathcal{O})$ is a ringed Segal topos (see \cite{tvv}). In the same
way, for any derived stack $F$, $\mathbb{A}^{1}\times F$ represents a 
ring object $\mathcal{O}_{F}$ in $\mathbf{dSt}(k)/F$, and the pair
$(\mathbf{dSt}(k)/F,\mathcal{O}_{F})$ is a ringed Segal topos. The derived Segal category 
of $(\mathbf{dSt}(k)/F,\mathcal{O}_{F})$ is denoted by $L(F,\mathcal{O}_{F})$ (see \cite{tvv} for
a precise definition). 

We also have
a sub-Segal category $L_{qcoh}(F)\subset L(F,\mathcal{O}_{F})$ of quasi-coherent complexes, which
can be described in the following way. We write $F$ as the colimit of 
affine derived schemes $F\simeq Colim X_{i}$, where $X_{i}=\mathbb{R}Spec\, A_{i}$. 
For any $i$, we can consider the commutative dg-algebra $N(A_{i})$ obtained by normalizing 
$A_{i}$, and thus its unbounded Segal category of modules $L(N(A_{i})-Mod)$. One possible definition is
$$L_{qcoh}(F)=Holim_{i} L(N(A_{i})-Mod),$$
where this homotopy limit is taken in the model category of Segal categories (see \cite{tova} where this is done
using dg-categories). Like in the underived situation, 
any morphism $f : F \longrightarrow F'$ induces an adjunction
$$f^{*} : L_{qcoh}(F') \longrightarrow L_ {qcoh}(F,\mathcal{O}_{F}) \qquad
L_{qcoh}(F') \longleftarrow L_ {qcoh}(F,\mathcal{O}_{F}) : f_{*}^{qcoh},$$
but again $f_{*}^{qcoh}$ is not the functor induced by the direct image on the level 
of all complexes of $\mathcal{O}$-modules.

When $F$ is a derived Artin stack then the Segal category 
$L_{qcoh}(F)$ has a natural t-structure. By definition, an object
$E \in L_{qcoh}(F)$ is of non positive amplitude if for any flat morphism $u : X=\mathbb{R}Spec\, A \longrightarrow F$
with $A\in sk-CAlg$, the corresponding object $u^{*}(E)\in 
L_{qcoh}(X)\simeq L(N(A)-Mod)$ is cohomolgically concentrated in non positive degrees (as a 
complex over $k$). The heart of this t-structure is denoted by $QCoh(F)$ is called the category of
quasi-coherent sheaves over $F$. The natural morphism $it_{0}(F) \longrightarrow F$ induces
by direct images an equivalence 
$$QCoh(F)\simeq QCoh(t_{0}(F)).$$
In particular we see that the two Segal categories $L_{qcoh}(F)$ 
and $L_{qcoh}(t_{0}F)$ are both endowed with a t-structure and have the same
heart, but are different in general. In this way, the derived enhancement $F$
of $t_{0}F$ can also be considered as a modification of the derived category 
$D_{qcoh}(t_{0}F)$, keeping the heart unchanged. 

Let us assume that we have a pull-back square of derived Artin stacks
$$\xymatrix{
F' \ar[r]^{q} \ar[d]_-{v} & G'\ar[d]^-{u} \\
F'\ar[r]_-{p} & G.}$$
Then there exists a base change natural transformation
$$\alpha : u^{*}\circ p_{*}^{qcoh} \Rightarrow q_{*}^{qcoh}\circ v^{*}.$$
The natural transformation $\alpha$ is an equivalence in many interesting examples, 
for instance when $F$, $G$, $F'$ and $G'$ are all quasi-compact derived schemes with an affine diagonal
(more generally when $F$, $G$, $F'$ and $G'$ are all strongly quasi-compact and
$p_{*}^{qcoh}$ and $q_{*}^{qcoh}$ are of finite t-amplitude). That the base change formula is
satisfied without any flatness assumptions on $u$ is an important feature of derived
algebraic geometry. 

\item \textbf{Tangent and cotangent complexes:}
For any derived Artin stack $F$ there exists an object $\mathbb{L}_{F}\in L_{qcoh}(F)$
called the \emph{cotangent complex of $F$}. It is characterized by the following 
universal property: for any $A \in sk-CAlg$, any morphism $x : X=\mathbb{R}Spec\, A \longrightarrow F$
and any simplicial $A$-module $M$, there exists a natural equivalence between the
homotopy fiber of $F(A\oplus M) \longrightarrow F(A)$ at $x$, and 
$Map_{N(A)-Mod}(x^{*}(\mathbb{L}_{F}),N(M))$ (here $N$ is the normalization functor
going from simplicial algebras and simplicial modules to dg-algebras and dg-modules). Of course, 
when $F=i(X)$ for $X$ a scheme, then $\mathbb{L}_{X} \in D_{qcoh}(X)$ is the usual
cotangent complex of $X$ (e.g. as defined in \cite{il}). For 
$x : X=\mathbb{R}Spec\, A \longrightarrow F$ we define \emph{the tangent complex of $F$ at $x$} to be
$$\mathbb{T}_{x}F:=\mathbb{R}\underline{Hom}(x^{*}(\mathbb{L}_{F}),A),$$
the dual of $x^{*}(\mathbb{L}_{F})$. 

For any morphism of derived Artin stacks $f : F \longrightarrow F'$ we define
a relative cotangent complex  by the following triangle in $L_{qcoh}(F)$
$$f^{*}(\mathbb{L}_{F'}) \longrightarrow \mathbb{L}_{F} \longrightarrow \mathbb{L}_{F/F'}.$$
It can be shown that $f$ is smooth if and only if it is locally of finite presentation and if
$\mathbb{L}_{F/F'}$ is of non negative t-amplitude (i.e. if $[\mathbb{L}_{F/F'},E[i]]=Ext^{i}(
\mathbb{L}_{F/F'},E)=0$ for all $i>0$ and for all $E$ belonging to the non positive part 
of the t-structure on $L_{qcoh}(F)$). In the same way $f$ is etale if and
only if it is locally of finite presentation and  $\mathbb{L}_{F/F'}\simeq 0$ (see \cite[2.2.5]{hagII}). 

\item \textbf{The virtual structure sheaf:} Let $F$ be a derived Artin stack and $t_{0}F$ its truncation. 
For a smooth morphism $U=\mathbb{R}Spec\, A \longrightarrow F$, we can consider the
graded $\pi_{0}(A)$-module $\pi_{*}(A)$ as a graded quasi-coherent sheaf on $Spec\, \pi_{0}(A)$. 
When $U$ varies over smooth morphisms to $F$, the various graded 
quasi-coherent sheaves $\pi_{*}(A)$ glue together and descend to a global graded 
quasi-coherent sheaf $\pi_{*}(\mathcal{O}^{virt})$ on the stack $t_{0}F$. This graded sheaf is
called the \emph{virtual structure sheaf of $F$}. In any case it is an important invariant 
living on $t_{0}F$ and remembering some information about the derived enhancement $F$ of $t_{0}F$. 

\end{enumerate}

\subsection{Some examples}

\begin{enumerate}

\item \textbf{Derived fiber products of schemes and stacks:} As we have said
the natural inclusion functor
$$i : \mathbf{St}(k) \longrightarrow \mathbf{dSt}(k)$$
is not left exact an in particular does not preserve fiber products. Therefore, 
a very first example of derived Artin stacks is given by considering a 
diagram of Artin stacks 
$$\xymatrix{
 & F \ar[d] \\
 G \ar[r] & H,}$$
and then considering $i(F)\times_{i(H)}i(G)$. The natural morphism 
$$i(F\times_{H}G) \longrightarrow i(F)\times_{i(H)}i(G)$$
is in general not an equivalence, thought the induced morphism on the truncations
$$t_{0}(i(F\times_{H}G)) \longrightarrow t_{0}(i(F)\times_{i(H)}i(G))$$
is an equivalence of stacks. Therefore, the derived Artin stack 
$i(F)\times_{i(H)}i(G)$ is a derived enhancement (in the sense of Def. \ref{d4-}) of the usual fiber products
of stacks. 

A very simple, but fundamental, example is when $F$, $G$ and $H$ are all affine schemes
given by a diagram of commutative k-algebras
$$\xymatrix{
 & A \\
 B & \ar[l] C\ar[u].}$$
Then, the derived stack $i(F)\times_{i(H)}i(G)$ is $\mathbb{R}Spec\, (A\otimes ^{\mathbb{L}}_{C}B)$, where
$A\otimes ^{\mathbb{L}}_{C}B$ is the derived tensor product computed in simplicial commutative rings. 
We see that $\pi_{i}(A\otimes ^{\mathbb{L}}_{C}B)=Tor_{i}^{C}(A,B)$, and thus that 
the virtual structure sheaf (see point $(6)$ of the last section)  on $Spec\, (A\otimes_{C}B)$ is $Tor_{*}^{C}(A,B)$. When
$H$ is a regular scheme scheme and $F$ and $G$ are closed subschemes intersecting properly, then 
this virtual structure sheaf on $Spec\, (A\otimes_{C}B)$ precisely compute the correct
intersection number. 

An interesting general construction involving fiber products is the inertia stack. Recall that 
for a stack $F$, the inertia stack is $I_{F}:=Map(S^{1},F)=F\times_{F\times F}F$, and that it
classifies objects endowed with an automorphism in $F$. Considering $F$ as 
a derived stack $i(F)$ we get a derived inertia stack
$$\mathbb{R}I_{F}:=I_{i(F)}:=i(F)\times_{i(F)\times i(F)}i(F) \in \mathbf{dSt}(k).$$
The derived stack $\mathbb{R}I_{F}$ is of course a derived enhancement of 
the stack $I_{F}$, and is naturally a derived group stack over $F$. When $F$ is a scheme,
then $I_{F}=F$ and thus $\mathbb{R}I_{F}$ is a natural non trivial derived enhancement
of $F$. For instance, when $F=Spec\, A$ is an affine scheme, then we have
$$\mathbb{R}I_{F}=\mathbb{R}Spec\, (A\otimes^{\mathbb{L}}_{A\otimes^{\mathbb{L}}A}A)$$
whereas $I_{F}=F$.  As $A\otimes^{\mathbb{L}}_{A\otimes^{\mathbb{L}}A}A$ is known to be
Hochschild homology of $A$, $\mathbb{R}I_{F}$ is some sort of 
global version of Hochschild homology of the stack $F$.

\item \textbf{Derived schemes vs dg-schemes:} Assume that $k$ has characteristic zero. 
A dg-scheme 
is a pair $(X,\mathcal{A}_{X})$, consisting of a scheme $X$ together with 
a sheaf of quasi-coherent commutative dg-$\mathcal{O}_{X}$-algebras $\mathcal{A}_{X}$
such that $\mathcal{A}^{0}_{X}=\mathcal{O}_{X}$ and $\mathcal{A}^{i}_{X}=0$ for $i>0$. This notion has been introduced some
years ago as models for derived schemes in order to construct derived moduli spaces
(see \cite{ck1,ck2}). There exists an obvious notion of morphisms of dg-schemes and
of a quasi-isomorphisms between them. The Segal category of dg-schemes, obtained
by localizing along quasi-isomorphisms will be denoted by $L(dg-Sch)$. As explained in 
\cite{tv3} it is possible to construct a functor
$$\Theta : L(dg-Sch) \longrightarrow \mathbf{dSt}(k),$$
and taking its values inside the sub-Segal category of derived schemes. 
When $X=Spec\, A$ is affine and $\mathcal{A}_{X}$ is given by 
a non positively graded commutative dg-algebra $A_{X}$, then 
$\Theta(X,\mathcal{A}_{X})$ is defined as being $\mathbb{R}Spec\, D(A_{X})$, where
$D(A_{X})$ denotes the commutative simplicial algebra obtained by denormalization from $A_{X}$. 
In general $\Theta(X,\mathcal{A}_{X})$ is defined using some covering 
of $X$ by affine schemes and gluing in a rather straightforward way. 

Essentially nothing is known about the functor $\Theta$ but I tend to think it is not well behaved (e.g. is not
fully faithful). The reason for this feeling is that by definition for any dg-scheme $(X,\mathcal{A}_{X})$
there exists a natural closed immersion of dg-schemes
$$(X,\mathcal{A}_{X}) \longrightarrow (X,\mathcal{O}_{X})=X.$$
Moreover, any morphism between dg-schemes $(X,\mathcal{A}_{X}) \longrightarrow (Y,\mathcal{A}_{Y})$
induces a commutative square of derived schemes
$$\xymatrix{
(X,\mathcal{A}_{X}) \ar[r] \ar[d] & X \ar[d] \\
(Y,\mathcal{A}_{Y}) \ar[r] & Y.}$$
This most probably implies that $\Theta$ is not essentially surjective on derived schemes, because 
there are no reasons for a given derived scheme $Z$ to be embeddable as a closed sub derived scheme
of some ambiant scheme (though such an embedding always exists locally). 
If we think of derived schemes as being somehow analogs of formal schemes, 
the ones that have such embedding are analogs of the algebraizable formal schemes. Moreover, 
following this analogy, the morphisms of derived schemes 
$\Theta(X,\mathcal{A}_{X}) \longrightarrow \Theta(Y,\mathcal{A}_{Y})$
compatible with a morphism $X \longrightarrow Y$ as above are analogs of the algebraizable 
morphisms between formal schemes. This analogy also suggests that 
$\Theta$ is probably not full. 

Even though the functor $\Theta$ is probably not well behaved it can be used to 
produce examples of derived schemes as already some interesting example of
dg-schemes have been constructed. This is for example the case of the derived
Quot and Hilbert schemes, defined in \cite{ck1,ck2}. These dg-schemes have been defined
in a rather ad-hoc manner, and they have not been constructed to represent
any moduli functors (actually, dg-schemes do not seem to be well suited for the functorial point of view,
due to the lack of a model structure on them that would allow to compute the
maps in the localization $L(dg-Sch)$). 
A natural question is therefore to describe 
moduli characterizations of derived schemes arising from  dg-schemes
via the functor $\Theta$. The two major examples are the derived schemes $\Theta(\mathbb{R}Quot(X))$ and
$\Theta(\mathbb{R}Hilb(X))$. The case of $\Theta(\mathbb{R}Quot(X))$ has been 
recently studied by J. Gorski in his thesis \cite{go}.

To conclude this comparison between dg-schemes and derived schemes: dg-schemes seem to be
only approximation of the correct notion of what a derived scheme is. Moreover, it seems there
is nothing doable with dg-schemes that can not be done with derived schemes, but there are
things doable with derived schemes that can not be done with dg-schemes (e.g. the construction of
some derived moduli such as vector bundles on some smooth proper scheme, or 
also having functorial description of these derived moduli). It is therefore reasonable
to suggest to simply forget about the notion of dg-schemes (and this has actually been suggested
once to me by M. Kontsevich). 

\item \textbf{Linear derived stacks:} As in example $(2)$ of \S 3.2 it is possible to define 
the notion of linear stacks
in the context of derived stacks. Let $F$ be any derived Artin stack, and 
$E\in L_{qcoh}(F)$. We define a derived stack $\mathbb{RV}(E)$ over $F$ by
$$\begin{array}{cccc}
\mathbb{RV}(E) : & \mathbf{dSt}(k)/F & \longrightarrow & Top \\
 & (f : F' \rightarrow F) & \mapsto & Map_{L(F',\mathcal{O}_{F'})}(f^{*}(E),\mathcal{O}_{F'}).
\end{array}$$
The derived stack $\mathbb{RV}(E)$ is 
called \emph{the linear stack associated to $E$}. By construction, it is characterized by the
following universal property
$$\pi_{0}(Map_{\mathbf{dSt}(k)/F}(F',\mathbb{RV}(E)))\simeq Ext^{0}(f^{*}(E),\mathcal{O}_{F'}),$$
for any $f : F' \longrightarrow F$ in $\mathbf{dSt}(k)/F$, and where the $Ext^{0}$ is computed in the
derived category of complexes of $\mathcal{O}_{F'}$-modules. 

The stack $\mathbb{RV}(E)$ is a derived Artin stack if 
$E$ is perfect (i.e. its pull-backs to any derived affine scheme $X=\mathbb{R}Spec\, A 
\longrightarrow F'$ is a compact object
in $D(N(A))$, the derived category of $N(A)$-dg-modules), and the morphism
$\mathbb{RV}(E) \longrightarrow F$ is then strongly of finite presentation.
The main difference with the notion of linear stacks in the underived situation is that 
$\mathbb{RV}(E)$ depends on the full complex $E$, and not only on its positive part. We 
have $t_{0}\mathbb{RV}(E)\simeq \mathbb{V}(E)$, and thus
$\mathbb{RV}(E)$ is a natural derived enhancement of $\mathbb{V}(E)$.

For any derived Artin stack $F$, locally of finite presentation, the
cotangent complex $\mathbb{L}_{F}$ is perfect. The derived Artin stack
$\mathbb{RV}(\mathbb{L}_{F})$ can then be identified with 
the derived  tangent stack $\mathbb{R}TF:=\mathbb{R}\mathbf{Map}(Spec\, k[\epsilon],F)$, defined
as the derived stack of morphisms from $Spec\, k[\epsilon]$ to $F$ (see point $(4)$ below). It is important to note that 
when $F$ is an Artin stack (i.e. $F\simeq it_{0}F$), then 
$\mathbb{R}TF$ is no longer an Artin stack except when $F$ is itself smooth. Therefore, 
even though $F$ is an Artin stack, $\mathbb{R}TF$ is in general a non trivial 
derived enhancement of the usual tangent stack $TF$.

\item \textbf{Derived mapping stacks:} As the Segal category $\mathbf{dSt}(k)$  is a Segal topos it
has internal Homs objects. These objects are denoted by 
$\mathbb{R}\mathbf{Map}$, in order to avoid confusions with the one defined 
for underived stacks. The functor $i : \mathbf{St}(k) \longrightarrow \mathbf{dSt}(k)$
does not commute with taking internal Homs, but its right adjoint
$t_{0}$ does. In particular, for two Artin stacks $F$ and $G$, the derived
stack $\mathbb{R}\mathbf{Map}(i(F),i(G))\in \mathbf{dSt}(k)$ is such that 
$t_{0}\mathbb{R}\mathbf{Map}(i(F),i(G))\simeq \mathbf{Map}(F,G)$, and is thus a derived
enhancement of $\mathbf{Map}(F,G)$. This provides a systematic way to construct non trivial examples
of derived stacks starting from underived stacks. 

For a stack $X$ and an Artin stack $F$ there exist criteria ensuring that 
$\mathbb{R}\mathbf{Map}(i(X),i(F))$ is a derived Artin stack. The most powerful 
follows from Lurie's representability criterion (see \cite{lu}), and states that this is the case
as soon as $X$ is a flat and proper scheme and $F$ is an Artin stack locally of finite presentation. 
A simpler, but less prowerful, criterion is given in \cite[App. C]{hagII}, and states that 
this is the case if it is already known that $ \mathbf{Map}(X,F)$ is an Artin stack and under some
additional mild conditions. These two criteria can be used to prove the existence of the following
derived Artin stacks:
\begin{enumerate}
\item For a finite and connected simplicial set $K$, the derived stack 
$$\mathbb{R}\mathbf{Loc}_{n}(K):=\mathbb{R}\mathbf{Map}(K,i(BGl_{n}))$$ 
is a derived Artin stack strongly of finite presentation and
is called the \emph{derived moduli stack of rank n local systems on $K$} (or on its geometric realization).
Its truncation 
$$\mathbf{Loc}_{n}(K):=\mathbf{Map}(K,BGl_{n})$$
 is the usual Artin 1-stack of 
rank n local systems on $K$, or equivalently of rank n linear representations 
of the group $\pi_{1}(K)$. Given a local system $V \in \mathbf{Loc}_{n}(K)(k)$, the tangent 
complex of $\mathbb{R}\mathbf{Loc}_{n}(K)$ at the point $V$ is 
$$\mathbb{T}_{V}\mathbb{R}\mathbf{Loc}_{n}(K)\simeq C^{*}(K,V\otimes V^{\vee})[1],$$
where $C^{*}(K,V\otimes V^{\vee})$ is the complex of cohomology of $K$ with 
coefficients in the local system of endomorphisms of $V$. This implies in particular that 
$\mathbb{R}\mathbf{Loc}_{n}(K)$ depends on more than $\pi_{1}(K)$ alone, and also
captures higher homotopical invariants of $K$. 
\item For a flat and proper scheme $X$, the derived stack
$$\mathbb{R}\mathbf{Vect}_{n}(X):=\mathbb{R}\mathbf{Map}(i(X),i(BGl_{n}))$$
is a derived Artin stack strongly of finite presentation and 
is called the \emph{derived moduli stack of rank n vector bundles on $X$}.
Its truncation 
$$\mathbf{Vect}_{n}(X):=\mathbf{Map}(X,BGl_{n})$$
 is the usual Artin 1-stack of 
rank n vector bundles on $X$. Moreover, for a vector bundle $V$ on $X$ we have
$$\mathbb{T}_{V}\mathbb{R}\mathbf{Vect}_{n}(X)\simeq C^{*}(X,V\otimes V^{\vee})[1].$$
\item Let $k$ be a field of characteristic zero and $X$ a smooth and projective scheme over $Spec\, k$. 
We consider $X_{DR}$ as defined in point $(4)$ of \S 3.3. The derived stack  
$$\mathbb{R}\mathbf{Vect}_{n}(X_{DR}):=\mathbb{R}\mathbf{Map}(i(X_{DR}),i(BGl_{n}))$$
is a derived Artin stack strongly of finite presentation and 
is called the \emph{derived moduli stack of rank n flat vector bundles on $X$}.
Its truncation 
$$\mathbf{Vect}_{n}(X_{DR}):=\mathbf{Map}(X_{DR},BGl_{n})$$
is the usual Artin 1-stack of 
rank n flat vector bundles on $X$. Moreover, for a flat vector bundle $V$ on $X$ we have
$$\mathbb{T}_{V}\mathbb{R}\mathbf{Vect}_{n}(X_{DR})\simeq C^{*}(X_{DR},V\otimes V^{\vee})[1],$$
where $C^{*}(X_{DR},V\otimes V^{\vee})$ is the complex of de Rham cohomology of $W$ with coefficients
in the flat bundle of endomorphisms of $V$. 

\item Let $\overline{\mathcal{M}}^{pre}_{g,n}$ be the 
Artin $1$-stack of prestable curves of 
genus $g$ and 
with $n$ marked points, and let $X$ be a smooth and proper scheme. We consider
the universal prestable curve $\overline{\mathcal{C}}_{g,n} \longrightarrow \overline{\mathcal{M}}^{pre}_{g,n}$. 
We define the derived stack of prestable maps to be
$$\mathbb{R}\overline{\mathcal{M}}^{pre}_{g,n}(X):=
\mathbb{R}\mathbf{Map}_{\mathbf{dSt}(k)/i(\overline{\mathcal{M}}^{pre}_{g,n})}
(i(\overline{\mathcal{C}}_{g,n}),X\times
i(\overline{\mathcal{M}}^{pre}_{g,n})),$$
where $\mathbb{R}\mathbf{Map}_{\mathbf{dSt}(k)/i(\overline{\mathcal{M}}^{pre}_{g,n})}$ denotes the
internal Homs of the comma Segal category of derived stacks over $i(\overline{\mathcal{M}}^{pre}_{g,n})$. 

The derived stack $\mathbb{R}\overline{\mathcal{M}}_{g,n}^{pre}(X)$
is naturally a derived stack over $i(\overline{\mathcal{M}}_{g,n}^{pre})$. Moreover, 
for any $Y=\mathbb{R}Spec\, B \longrightarrow i(\overline{\mathcal{M}}^{pre}_{g,n})$ we have
$$\mathbb{R}\overline{\mathcal{M}}^{pre}_{g,n}(X)\times_{i(\overline{\mathcal{M}}^{pre}_{g,n})}Y\simeq 
\mathbb{R}\mathbf{Map}_{\mathbf{dSt}(k)/Y}(i(\overline{\mathcal{C}}_{g,n})\times
_{i(\overline{\mathcal{M}}^{pre}_{g,n})}Y,X\times Y).$$
This implies that the morphism $\mathbb{R}\overline{\mathcal{M}}^{pre}_{g,n}(X)\longrightarrow
i(\overline{\mathcal{M}}_{g,n}^{pre})$ is a relative derived Artin stack and thus that 
$\mathbb{R}\overline{\mathcal{M}}^{pre}_{g,n}(X)$ is a derived Artin stack. 
The derived
stack of stable maps $\mathbb{R}\overline{\mathcal{M}}_{g,n}(X)$ 
is the open derived substack of $\mathbb{R}\overline{\mathcal{M}}_{g,n}^{pre}(X)$ consisting
of stable maps. In other words, there exists a cartesian square
$$\xymatrix{
i(\overline{\mathcal{M}}_{g,n}(X)) \ar[r] \ar[d] & \mathbb{R}\overline{\mathcal{M}}_{g,n}(X) \ar[d] \\
i(\overline{\mathcal{M}}_{g,n}^{pre}(X)) \ar[r] & \mathbb{R}\overline{\mathcal{M}}_{g,n}^{pre}(X),}$$
where $\overline{\mathcal{M}}_{g,n}(X) \subset \overline{\mathcal{M}}_{g,n}^{pre}(X)$ 
is the substack of stable maps. 

The truncation of $\mathbb{R}\overline{\mathcal{M}}_{g,n}(X)$ is by construction the usual 
stack of stable maps to $X$, and therefore we see that  
$\mathbb{R}\overline{\mathcal{M}}_{g,n}(X)$ is a derived Deligne-Mumford stack. 
The tangent of $\mathbb{R}\overline{\mathcal{M}}_{g,n}(X)$ at a 
morphism $f : C \longrightarrow X$, is given by 
$$\mathbb{T}_{f}\mathbb{R}\overline{\mathcal{M}}_{g,n}(X)\simeq 
C^{*}(C,T_{C}(-\sum x_{i}) \rightarrow f^{*}(T_{X}))[1],$$
where $T_{C}(-\sum x_{i}) \rightarrow f^{*}(T_{X})$ is a complex of sheaves on $C$ concentrated in degrees
$[0,1]$. 
As a consequence of point $(6)$ of \S 4.3 we immediately get 
a virtual structure sheaf $\pi_{*}(\mathcal{O}^{virt})$ on 
the usual stack of prestable maps $\overline{\mathcal{M}}_{g,n}^{pre}(X)$ and therefore
on the usual stack of stable maps $\overline{\mathcal{M}}_{g,n}(X)$.
\end{enumerate}

\item \textbf{Objects in a dg-category:} For a saturated dg-category $T$, 
the locally geometric Artin  stack $\mathbf{Parf}_{T}$ described in example $(4)$ of \S 3.2 has a natural 
derived enhancement denoted by $\mathcal{M}_{T}$ (see \cite{tova} for a precise definition
of the derived stack $\mathcal{M}_{T}$). It can also be proved that 
$\mathcal{M}_{T}$ is locally geometric (i.e. union of open derived Artin sub-stacks locally of 
finite presentation). In fact, the original proof of the local 
geometricity of $\mathbf{Parf}_{T}$ is deduced from the one of $\mathcal{M}_{T}$ which 
is somehow simpler, has explicit computations of cotangent complexes help 
proving the existence of a smooth atlas. For a given object $E$ in $T$, 
the tangent complex is given by
$$\mathbb{T}_{E}\mathcal{M}_{T}\simeq Ext^{*}(E,E)[1],$$
where $Ext^{*}(E,E)=T(E,E)$ is the complex of endomorphisms of $E$. 

An important consequence of the local geometricity of $\mathcal{M}_{T}$ is the
existence of a local geometric derived stack $\mathbb{R}\mathbf{Parf}(X)$, of 
perfect complexes on some smooth and proper scheme. The derived stack 
$\mathbb{R}\mathbf{Parf}(X)$ is of course a derived enhancement of $\mathbf{Parf}(X)$
described in example $(4)$ of \S 3.2. Inside the stack $\mathbf{Parf}(X)$ sits as an open sub-stack 
$\mathbf{Coh}(X)$ the 1-stack of coherent sheaves on $X$. As 
the stack $\mathbf{Parf}(X)$ and $\mathbb{R}\mathbf{Parf}(X)$ have the same
topology (and in particular the same open substacks), there exists a unique
open derived sub-stack $\mathbb{R}\mathbf{Coh}(X)\subset \mathbb{R}\mathbf{Parf}(X)$ such that the
diagram 
$$\xymatrix{
i(\mathbf{Coh}(X)) \ar[r] \ar[d] & i(\mathbf{Parf}(X)) \ar[d] \\
\mathbb{R}\mathbf{Coh}(X) \ar[r] &  \mathbb{R}\mathbf{Parf}(X)}$$
is cartesian. 

This implies that $\mathbb{R}\mathbf{Coh}(X)$ is a derived Artin stack. 
The derived stack $\mathbb{R}\mathbf{Coh}(X)$ itself contains
a derived open sub-stack $\mathbb{R}\mathbf{Vect}_{n}(X)$ of vector bundles on $X$ of rank $n$.
This provides another direct proof of the geometricity of $\mathbb{R}\mathbf{Vect}_{n}(X)$
without refering to any representability criterion. An interesting example
is $\mathbb{R}\mathbf{Pic}(X)=\mathbb{R}\mathbf{Vect}_{1}(X)$, the derived Picard stack 
of $X$. Indeed, the truncation $t_{0}\mathbb{R}\mathbf{Pic}(X)$ is the usual Picard
stack of line bundles on $X$, and thus is smooth. However, though its
truncation is smooth it is not true that 
$\mathbb{R}\mathbf{Pic}(X)=it_{0}(\mathbb{R}\mathbf{Pic}(X))$ as this can be seen on the
tangent complexes. This example shows that the usual intuition that moduli spaces
are singular because of the existence of a non trivial derived structure is not always
true in practice.

\item \textbf{Dg-categories:} The stack $\mathbf{dgCat}^{sat}$ has a natural derived enhancement
$\mathbb{R}\mathbf{dgCat}^{sat}$ defined in the following way. For any $A \in sk-CAlg$, 
we consider the commutative dg-algebra $N(A)$ obtained by normalizing $A$. 
The category of $N(A)$-dg-modules has a natural symmetric monoidal structure, and therefore
it make sense to talk about $N(A)$-dg-categories. Moreover, the notion of 
being saturated naturally extends from dg-categories over $k$ to dg-categories over $N(A)$. 
The functor sending $A$ to the nerve of the category of quasi-equivalences between 
saturated $N(A)$-dg-categories is denoted by $\mathbb{R}\mathbf{dgCat}^{sat}$. 
We clearly have $t_{0}\mathbb{R}\mathbf{dgCat}^{sat}\simeq \mathbf{dgCat}^{sat}$. 
As in the underived case, I believe that $\mathbb{R}\mathbf{dgCat}^{sat}$ is a locally geometric
derived stack, and I think a direct approach using Lurie's representability criterion should
be possible. 

\begin{Q}
Is the derived stack $\mathbb{R}\mathbf{dgCat}^{sat}$ locally geometric ?
\end{Q}

Of course a positive answer to this question would also provide a positive answer to question \ref{Q1}. 
Naturally, it is expected that for a given saturated dg-category $T$ the tangent complex is given by
$$\mathbb{T}_{T}\mathbb{R}\mathbf{dgCat}^{sat}\simeq HH(T)[2],$$
where $HH(T)$ is the full Hochschild cohomology of $T$. 

\end{enumerate}

\subsection{Some developments}

\begin{enumerate}

\item \textbf{Representability criterion:} Probably the most important recent development in the
theory of derived stack is the representability criterion proved by J. Lurie in \cite{lu}, that we already mentioned
several times but that we will not reproduce here. 
It is a generalization of the standard Artin's representability criterion for algebraic spaces and
1-stacks. However, the criterion in the derived setting is simpler as the part concerning 
having a good infinitesimal theory is now truly a property of the moduli functor and not
an extra structure as we explained during the introduction of this section (see \S 4.1). This criterion is
extremely powerful, though it is not always very easy to check the infinitesimal
properties in practice, and it is sometimes easier to prove directly the
geometricity and then deduced the infinitesimal theory from it (this is what is done
for example in several examples in \cite{hagII}). 

\item \textbf{Formal theory and derived inertia stacks:} Assume that $k$ has characteristic zero. Let $F$ be 
a derived Artin stack locally of finite presentation and $\mathbb{R}I_{F}=\mathbb{R}\mathbf{Map}(S^{1},F)$ be its derived
inertia stack. The composition of loops makes $\mathbb{R}I_{F}$ into a group object over
$F$. In particular for any point $x : Spec\, K \longrightarrow F$, with $K$ a field, 
we obtain $\mathbb{R}I_{F}\times_{F}Spec\, K$, which is a derived group Artin stack over $Spec\, K$. 
This group object has a tangent Lie algebra $L_{x}$ (which is well defined 
in the homotopy category of dg-Lie algebras over $K$). The precise relation between derived group stack
and dg-Lie algebra has not been investigated yet, and there might be some foundational 
work to be done to explain what $L_{x}$ truly is. In any case, I will assume that we know how to do 
this. It is easy to see that, as a complex, $L_{x}$ is naturally quasi-isomorphic
to $\mathbb{T}_{x}F[-1]$, the shifted tangent complex of $F$ at $x$. Therefore, 
we obtain a natural structure of dg-Lie algebra (or at least $L_{\infty}$-Lie algebra) on $\mathbb{T}_{x}F[-1]$.
From this dg-Lie algebra $L_{x}$ we can define a formal derived moduli functor, defined
on the  category of augmented Artinian dg-algebras over $K$ (see \cite{hin})
$$MC(L_{x}) : dg-Art/K \longrightarrow SSet.$$
On the other hand, the restriction of the derived stack $F$ on $dg-Art/K$, pointed at $x$ also 
provides a functor
$$\widehat{F}_{x} : dg-Art/K \longrightarrow SSet$$
which is by definition the formal completion of $F$ at the point $x$. It is expected that 
the two formal derived stacks $MC(L_{x})$ and $\widehat{F}_{x}$ are in fact equivalent. 
In other words, the tangent complex $\mathbb{T}_{x}F$ together with the
dg-Lie algebra structure on $\mathbb{T}_{x}F[-1]$ determines the formal completion of the derived 
stack at $x$. This statement seems to have been 
proved for derived schemes and derived Deligne-Mumford stacks
as it can be essentially reduced to the case of an affine derived stack which is somehow treated in 
\cite{hin,man}. However, the general statement for higher derived stacks does not seem to be 
known. 

\begin{Q}
Compare the two formal derived stacks $MC(L_{x})$ and $\widehat{F}_{x}$. 
\end{Q}

The local picture around the point $x$ also has a global counter-part, as the
group object $\mathbb{R}I_{F}$ has a global sheaf of quasi-coherent dg-Lie algebras
$L$ on $F$, whose underlying complex is $\mathbb{T}_{F}[-1]$
the shifted tangent complex of $F$. The sheaf of dg-Lie algebras $L$ on $F$
is supposed to control the formal completion of $F\times F$ along the diagonal (thought the
precise meaning of this in the stack context is not completely clear).

\item \textbf{Virtual fundamental classes:}  Let $F$ be a derived Artin stack, $t_{0}F$ its truncation
and $\pi_{*}(\mathcal{O}_{F}^{virt})$ the graded virtual structure sheaf on $t_{0}F$ as
defined in point $(6)$ of \S 4.2. When the quasi-coherent sheaves $\pi_{i}(\mathcal{O}_{F}^{virt})$
are all coherent and vanish for $i$ big enough, we can define a virtual 
class in $G$-theory
$$[\mathcal{O}_{F}]^{virt}:=\sum(-1)^{i}[\pi_{i}(\mathcal{O}_{F}^{virt})]\in G_{0}(t_{0}F).$$
Note that when $it_{0}F\simeq F$ then $[\mathcal{O}_{F}]^{virt}=[\mathcal{O}_{F}]$. 

The condition that the sheaves $\pi_{i}(\mathcal{O}_{F}^{virt})$
are coherent and vanish for $i$ big enough is not often satisfied and is a rather strong condition. 
It is known to be satisfied when $k$ is noetherian and of characteristic zero, $F$ is locally of finite presentation and
the cotangent complex $\mathbb{L}_{F}$ is of amplitude contained in $[-1,\infty[$. When $F$ satisfies these
two conditions we will say that $F$ is \emph{quasi-smooth}. 

Assume now that $k$ is noetherian of characteristic zero and
and that $F$ is a quasi-smooth derived Deligne-Mumford stack. On one hand we have 
the virtual class in $G$-theory $[\mathcal{O}_{F}]^{virt}\in G_{0}(t_{0}F)$, from which 
we can construct via the Grothendieck-Riemann-Roch transformation a class
in the rational Chow groups $\tau([\mathcal{O}_{F}]^{virt})\in CH_{*}(t_{0}F)_{\mathbb{Q}}$. 
On the other hand, we can pull-back the cotangent complex
via the morphism $j : it_{0}F \longrightarrow F$. The complex
$j^{*}(\mathbb{L}_{F})$ is a perfect obstruction theory of amplitude $[-1,0]$ in the sense
of \cite{bf}, and thus we can also construct a virtual fundamental class 
$[F]^{virt}\in CH_{*}(t_{0}F)_{\mathbb{Q}}$. As far as I know the following question is still open.

\begin{Q}\label{Q3}
What is the relation between $\tau([\mathcal{O}_{F}]^{virt})$ and $[F]^{virt}$ ?
\end{Q}

It seems that it is expected that these two classes only differ by a $Td(T_{F}^{virt})$, where
$T_{F}$ is the virtual tangent sheaf defined to be the dual of 
$j^{*}(\mathbb{L}_{F})$  (see \cite{bf}). Some results in that direction are proved in \cite{jo}.

Finally, when $F$ is not quasi-smooth anymore, and more generally when 
the virtual structure sheaf has infinite non zero sheaves, it is very much unclear 
how to use this virtual sheaf in order to get interesting invariants generalizing the
virtual class.  

\item \textbf{A holomorphic Casson invariant:}
Suppose now that $k=\mathbb{C}$, and that $X$ is a Calabi-Yau 3-fold. 
We consider $\mathbb{R}\mathbf{Coh}(X)$ and its derived open sub-stack
$\mathbb{R}\mathbf{Coh}^{st,\nu}(X)$ consisting of stable coherent sheaves
with some fixed numerical invariants $\nu \in K_{0}^{num}(X)$. We also assume that 
$\nu$ is chosen in such a way that semi-stable implies stable. Finally, we set 
$$\mathcal{M}(X):=[\mathbb{R}\mathbf{Coh}^{st,\nu}(X)/\mathbb{R}\mathbf{Pic}^{0}(X)],$$
the quotient derived stack of $\mathbb{R}\mathbf{Coh}^{st,\nu}(X)$ by the natural action 
of $\mathbb{R}\mathbf{Pic}^{0}(X)$ the derived group stack of line bundles of degree zero. 
The derived stack $\mathcal{M}(X,\nu)$ is a proper derived algebraic space. 
Moreover, the tangent complex at a coherent sheaf $E$ on $X$ can be seen to fit
in a triangle
$$C^{*}(X,\mathcal{O})[1] \longrightarrow C^{*}(X,E\otimes E^{\vee})[1] \longrightarrow 
\mathbb{T}_{E}\mathcal{M}(X,\nu).$$
Using the trace morphism $tr : C^{*}(X,E\otimes E^{\vee}) \longrightarrow C^{*}(X,\mathcal{O})$
we see that the above triangle splitts, and that the tangent complex at $E$ is given by
$$\mathbb{T}_{E}\mathcal{M}(X,\nu)\simeq C^{*}(X,E\otimes E^{\vee})^{0}[1],$$
where $C^{*}(X,E\otimes E^{\vee})^{0}$ is the kernel of the morphism $tr$. 

The conclusion is that $\mathcal{M}(X,\nu)$ is a quasi-smooth and proper derived algebraic space
over $Spec\, \mathbb{C}$, which is furthermore of virtual dimension zero. We can therefore
define a holomorphic Casson invariant to be the lenght of $\mathcal{M}(X,\nu)$ by the formula
$$\chi(X,\nu):=\sum(-1)^{i}h^{i}(\mathcal{M}(X,\nu),\mathcal{O}_{\mathcal{M}(X,\nu)}).$$
As the derived stack $\mathcal{M}(X,\nu)$ is proper and quasi-smooth, 
we see that $\chi(X,\nu)$ is well defined. It is also equal to 
$$p_{*}([\mathcal{O}_{\mathcal{M}(X,\nu)}]^{virt})\in G_{0}(Spec\, \mathbb{C})=\mathbb{Z},$$
where $p : t_{0}\mathcal{M}(X,\nu) \longrightarrow Spec\, \mathbb{C}$ is the projection and
$[\mathcal{O}_{\mathcal{M}(X,\nu)}]^{virt}$ is the virtual class in $G$-theory defined above in point 
$(6)$ of \S 4.2.

By construction, the invariant $\chi(X,\nu)$ counts the virtual number of
stable sheaves with numerical invariants $\nu$, with fixed determinants. It is natural to call
it the \emph{holomorphic Casson invariant}. It is probably very close to the one
defined in \cite{th}, as it surely satisfies the same deformation invariance property (this
is an application of the base change formula, point $(4)$ of \S 4.2 ). However, a precise 
comparison between these
two invariants requires an answer to the question \ref{Q3}. 

\item \textbf{Concerning the geometric Langlands correspondence:} As far as I understand
the geometric version of the Langlands correspondence predicts that for any 
smooth and projective curve $C$ over $k=\mathbb{C}$, the existence of 
an equivalence of triangulated categories
$$D(\mathbf{Vect}_{n}(C),\mathcal{D})\simeq D_{coh}(\mathbf{Loc}_{n}^{DR}(C))$$
where the left hand side is the derived category of $\mathcal{D}$-modules on 
the stack $\mathbf{Vect}_{n}(C)$ (with some finiteness conditions like being regular holonomic), and 
$\mathbf{Loc}_{n}^{DR}(C):=\mathbf{Map}(C_{DR},BGl_{n})$ is the stack 
of rank n flat bundles on $C$. I have recently learned from 
V. Lafforgue (and apparently this is a folklore knowledge shared by the
experts) that in order for this equivalence to have a chance to exist
the right hand side should rather be $D_{coh}(\mathbb{R}\mathbf{Loc}_{n}^{DR}(C))$, where 
$\mathbb{R}\mathbf{Loc}_{n}^{DR}(C)$ is the derived stack of rank n flat bundles 
discussed in example $(4-c)$ of \S 4.3. A striking example showing why this is so is when $n=1$. 

The stack $\mathbf{Vect}_{1}(C)$ is equivalent to $Pic^{0}(C) \times \mathbb{Z}\times K(\mathbb{G}_{m},1)$. 
On the other hand, the stack $\mathbf{Loc}_{1}^{DR}(C)$ is equivalent to 
$Pic^{0}(C)^{\dag}\times K(\mathbb{G}_{m},1)$, where $Pic^{0}(C)^{\dag}$ is the universal 
extension of the Jacobian $Pic^{0}(C)$ by the vector space $H^{0}(C,\Omega_{C}^{1})$. 
It is known that there exists an equivalence of triangulated categories (see \cite{la} for the first of these two equivalences)
$$D(Pic^{0}(C),\mathcal{D})\simeq D_{coh}(Pic^{0}(C)^{\dag}) \qquad 
D(\mathbf{Z},\mathcal{D})\simeq D_{coh}(K(\mathbb{G}_{m},1)).$$
Combining these two this shows that there exists an equivalence 
$$D(Pic^{0}(C)\times \mathbf{Z},\mathcal{D})\simeq D_{coh}(\mathbf{Loc}_{1}^{DR}(C)).$$
Therefore, we see that the part $D(K(\mathbb{G}_{m},1),\mathcal{D})$ is not reflected in 
$D_{coh}(\mathbf{Loc}_{1}^{DR}(C))$ and that the originial predicted equivalence does not seem to exist. 

Let $\mathbb{R}\mathbf{Loc}_{1}^{DR}(C)$ be the derived moduli stack of rank 1 flat bundles on $C$. 
It can be seen that we have
$$\mathbb{R}\mathbf{Loc}_{1}^{DR}(C)\simeq Pic^{0}(C)^{\dag}\times K(\mathbb{G}_{m},1) \times \mathbb{R}Spec\, \mathbb{C}[\mathbb{C}[1]],$$
where $\mathbb{C}[\mathbb{C}[1]]$ is the trivial square zero extension of $\mathbb{C}$ by 
$\mathbb{C}[1]$ (as a commutative dg-algebra it is freely generated by an element in degree -1). Moreover, it can be shown that 
there exists an equivalence of derived categories
$$D(K(\mathbb{G}_{m},1),\mathcal{D})\simeq D_{coh}(\mathbb{R}Spec\, \mathbb{C}[\mathbb{C}[1]]).$$
Indeed, $D(K(\mathbb{G}_{m},1),\mathcal{D})$ is equivalent to the derived category 
of $S^{1}$-equivariant complexes of $\mathbb{C}$-vector spaces, which is well known (via some bar-cobar construction) 
to be equivalent to the derived category of $\mathbb{C}[\mathbb{C}[1]]$-dg-modules. 

The conclusion is that the statement of the geometric Langlands correspondence is truly about 
the derived category of the derived stack $\mathbb{R}\mathbf{Loc}_{n}^{DR}(C)$. In the example above we also 
see that the "stacky part" $K(\mathbb{G}_{m},1)$ of $\mathbf{Vect}_{1}(C)$ correspond through the Langlands
correspondence to the "derived part" $\mathbb{R}Spec\, \mathbb{C}[\mathbb{C}[1]]$ of 
$\mathbb{R}\mathbf{Loc}_{1}^{DR}(C)$. This seems to be a general phenomenon, and explains somehow that
there exists some kind of duality between the stacky direction and the derived direction. On the infinitesimal level 
these two directions can be respectively observed as the
negative part and the  positive part of the tangent complex. In general, this duality between the stacky and the 
derived part can be understood in terms of characteristic cycles of $\mathcal{D}$-modules on higher stacks. 
Indeed, for a given Artin stack $F$, and a $\mathcal{D}$-module $E$ on $F$, the characteristic cylce of $E$
is supposed to live on the total cotangent stack of $F$. A reasonable candidate to be the cotangent stack
would be $\mathbb{V}(\mathbb{T}_{F})$, the linear stack associated to the tangent complex of $F$. 
But, when $F$ has a non trivial stacky direction (i.e. when it is at least a 1-Artin stack), then 
the complex $\mathbb{T}_{F}$ has non trivial negative cohomology sheaves, and thus
we have seen in example $(3)$ of \S 4.3 that $\mathbb{V}(\mathbb{T}_{F})$ has a natural non trivial derived enhancement
$\mathbb{RV}(\mathbb{T}_{F})$. The correct cotangent stack of $F$ is therefore a derived Artin stack, and the 
characteristic cycle of $E$ is now expected to live on $\mathbb{RV}(\mathbb{T}_{F})$.

\item \textbf{Categorified quantum cohomology:} Let $\mathbb{R}\overline{\mathcal{M}}_{g,n+1}(X)$ be the
derived stack of stable maps to a smooth and projective complex variety $X$ (see example $(4-d)$ of \S 4.3). Let us fix
a class $\beta \in H^{2}(X,\mathbb{Z})$, and 
let $\mathbb{R}\overline{\mathcal{M}}_{g,n+1}(X,\beta)$ the derived sub-stack of maps having 
$\beta$ as fundamental class. Forgetting the map to $X$ and eveluating at the marked points provide 
a natural diagram of derived stacks
$$\xymatrix{
\mathbb{R}\overline{\mathcal{M}}_{g,n+1}(X,\beta)\ar[r] \ar[d] & X \\
\overline{\mathcal{M}}_{g,n}\times X^{n}.}$$
This diagram induces by pull-back and push-forward a functor on the Segal categories of quasi-coherent complexes
$$L_{qcoh}(\overline{\mathcal{M}}_{g,n})\times L_{qcoh}(X)^{n} \longrightarrow L_{qcoh}(X).$$
This should be thought of as some kind of action of the system of Segal derived categories
$\{L_{qcoh}(\overline{\mathcal{M}}_{g,n})\}_{n,g}$ on $L_{qcoh}(X)$. The precise meaning of this
action must be made precise, and should be somehow an "action" of some kind of operad objects
in Segal categories (here it is preferable to use dg-categories instead of Segal categories in order to 
keep track of the linear structure). Note that the fact that this "action" satisfies the 
expected associativity axioms will follow from the base change formula (point $(4)$ of \S 4.2), showing the
importance to use the derived stacks $\mathbb{R}\overline{\mathcal{M}}_{g,n+1}(X,\beta)$
in the construction. 

This "action" of  $\{L_{qcoh}(\overline{\mathcal{M}}_{g,n})\}_{n,g}$ on $L_{qcoh}(X)$ 
can be thought of as a categorified version of quantum cohomolgy,
as passing from Segal categories to their Hochschild homology group would give back something
close to the quantum cohomology of $X$. What seems interesting with this construction is that 
the action of $\{L_{qcoh}(\overline{\mathcal{M}}_{g,n})\}_{n,g}$ on $L_{qcoh}(X)$ makes
sense even though $X$ is not smooth (one problem though is that the action of
$\{L_{qcoh}(\overline{\mathcal{M}}_{g,n})\}_{n,g}$ on $L_{qcoh}(X)$ does not preserve
bounded coherent complexes anymore). 

\end{enumerate}


\begin{thebibliography}{99}

\bibitem[An]{anel} M. Anel, in preparation. 

\bibitem[Be]{be} K. Behrend \textit{Derived l-adic categories for algebraic stacks}, 
Mem. Amer. Math. Soc. $\mathbf{163}$ (2003), no. 774, viii+93 pp. 

\bibitem[Be-Fa]{bf} K. Behrend, B. Fantechi, \textit{The intrinsic normal cone}, Invent. Math. $\mathbf{128}$
(1997), No. 1, 45-88. 

\bibitem[Ber1]{ber} J. Bergner, \textit{Three models for the homotopy theory of homotopy theories}, 
Preprint math.AT/0504334. 

\bibitem[Ber2]{ber2} J. Bergner, \textit{A characterization of fibrant Segal categories},
Preprint math.AT/0603400.  

\bibitem[Bo-VdB]{bovdb} A. Bondal, M. Van Den Bergh, \textit{Generators and representability of
functors in commutative and non-commutative geometry}, Mosc. Math. J. \textbf{3}
(2003), no. 1, 1--36.

\bibitem[Br]{bre} L. Breen, \textit{On classification of 2-gerbes and 2-stacks}, 
Ast\'erisque $\mathbf{225}$, Soc. Math. de France 1994.


\bibitem[Ci-Ka$1$]{ck1} I. Ciocan-Fontanine, M. Kapranov,
\textit{Derived Quot schemes}, Ann. Sci. Ecole Norm. Sup. (4) \textbf{34} (2001), 403-440.

\bibitem[Ci-Ka$2$]{ck2} I. Ciocan-Fontanine, M. Kapranov,
\textit{Derived Hilbert Schemes}, J. Amer. Math. Soc.  $\mathbf{15}$  (2002),  no. 4, 787-815.

\bibitem[Cis]{der} D.-C. Cisinski, \textit{Images directes cohomologiques dans les cat\'egories de mod\`eles}, 
Ann. Math. Blaise Pascal $\mathbf{10}$ (2003), no. 2, 195-244.

\bibitem[De-Xi]{de} B. Deng, J. Xiao, \textit{On Ringel-Hall algebras},
in Representations of finite dimensional algebras and related topics in Lie theory and geometry,
319--348, Fields Inst. Commun., $\mathbf{40}$, Amer. Math. Soc., Providence, RI, 2004.

\bibitem[DHI]{dhi} D. Dugger, S. Hollander, D. Isaksen, \textit{Hypercovers and simplicial presheaves}, Math. Proc. Cambridge Philos. Soc.
\textbf{136} (2004), 9-51.  

\bibitem[Dw-Ka1]{dk1} W. Dwyer,  D. Kan, \textit{Simplicial localization of categories},
J. Pure and Appl. Algebra $\mathbf{17}$ (1980), 267-284.

\bibitem[Dw-Ka2]{dk3} W. Dwyer,  D. Kan, \textit{Function complexes in homotopical algebra},
Topology $\mathbf{19}$ $(1980)$, $427-440$.

\bibitem[Go]{go} J. Gorski, \textit{Representability of the derived Quot functor}, thesis. 

\bibitem[Gr]{gr} A. Grothendieck, \textit{Pursuing stacks}, unpublished manuscript. 


\bibitem[Hin]{hin} V. Hinich, \textit{Formal stacks as dg-coalgebras}, J. Pure Appl. Algebra $\mathbf{162}$
(2001), No. 2-3, 209-250.


\bibitem[H-S]{hs} A. Hirschowitz, C. Simpson, \textit{Descente pour les $n$-champs}, 
preprint math.AG/$9807049$. 

\bibitem[Hol]{hol} S. Hollander, \textit{A homotopy theory for stacks},   
preprint math.AT/0110247.  
 
\bibitem[Ho1]{ho} M. Hovey, \textit{Model categories}, Mathematical surveys and monographs, Vol. $\mathbf{63}$, 
Amer. Math. Soc. Providence $1998$. 

\bibitem[Ho2]{ho2} M. Hovey, \textit{Model category structures on chain complexes of sheaves},
Trans. Amer. Math. Soc. \textbf{353} (2001), no. 6, 2441--2457.

\bibitem[Il]{il} L. Illusie, \textit{Complexe cotangent et d\'eformations I}, 
Lecture Notes in Mathematics \textbf{239}, Springer Verlag, Berlin, 1971. 

\bibitem[Ja1]{ja1} J. F. Jardine, \textit{Simplicial presheaves}, J. Pure and Appl. Algebra $\mathbf{47}$ (1987),  
35-87.  
  
\bibitem[Ja2]{ja2} J. F. Jardine, \textit{Stacks and the homotopy theory of simplicial sheaves},  
in \textit{Equivariant stable homotopy theory and related areas} (Stanford, CA, 2000).  
Homology Homotopy Appl. $\mathbf{3}$ (2001), No. 2, 361-384.   
  
\bibitem[Jo]{jo} R. Joshua, \textit{Riemann-Roch for algebraic stacks III: 
Virtual structure sheaves and virtual fundamental classes}, 
preprint available at http://www.math.ohio-state.edu/~joshua/pub.html.

\bibitem[Joy]{joy} A. Joyal, Letter to Grothendieck.  

\bibitem[Ka-Pa-Si]{kps} L. Katzarkov, T. Pantev, C. Simpson, \textit{Non-abelian mixed Hodge structures},
preprint math.AG/0006213.

\bibitem[Ka-Pa-To]{kpt} L. Katzarkov, T. Pantev, B. To\"en, \textit{Schematic homotopy types and
non-abelian Hodge theory},
preprint math.AG/0107129.

\bibitem[La]{la} G. Laumon, \textit{Transformation de Fourier generalis\'ee}, 
Preprint alg-geom/9603004.

\bibitem[La-Mo]{lm} G. Laumon and L. Moret-Bailly,  
\textit{Champs alg\'ebriques}, A series of Modern Surveys in 
Mathematics vol. $\mathbf{39}$, Springer-Verlag 2000.  

\bibitem[La-Ol]{laol} Y. Laszlo, M. Olsson, 
\textit{The six operations for sheaves on Artin stacks I: Finite Coefficients}, 
Preprint math.AG/0512097. 

\bibitem[Le1]{le} T. Leinster, \textit{A survey of definitions of n-category}, 
Theory and Applications of Categories 10 (2002), no. 1, 1-70.

\bibitem[Le2]{le2} T. Leinster, \textit{Higher Operads, Higher Categories}, 
London Mathematical Society Lecture Note Series $\mathbf{298}$, Cambridge University Press.

\bibitem[Lie]{lie} M. Lieblich, \textit{Moduli of complexes on a proper morphism}, 
to appear in Journal of Algebraic Geometry.

\bibitem[Lo-VdB1]{lvdb}  W. T. Lowen, M. 
Van den Bergh, \textit{Deformation theory of abelian categories}, 
to appear in Trans. Amer. Math. Soc.

\bibitem[Lo-VdB2]{lvdb2}  W. T. Lowen, 
M. Van den Bergh,\textit{Hochschild cohomology of abelian categories and ringed spaces}, 
Adv. in Math. $\mathbf{198}$ (2005), no. 1, 172-222.

\bibitem[Lu1]{lu} J. Lurie, \textit{Derived algebraic geometry}, thesis, available
at the author's home page. 

\bibitem[Lu2]{lu2} J. Lurie, \textit{On $\infty$-topoi}, Preprint math.CT/0306109.

\bibitem[Man]{man} M. Manetti, \textit{Extended deformation functors}, Int. Math. Res. Not. $\mathbf{14}$
 (2002), 719-756.

\bibitem[Ol1]{ol1} M. Olsson, \textit{F-isocrystals and homotopy types},
Preprint available at http://www.ma.utexas.edu/~molsson/.

\bibitem[Ol2]{ol2} M. Olsson, 
\textit{Towards non--abelian $P$--adic Hodge theory in the good reduction case}, Preprint
available at http://www.ma.utexas.edu/~molsson/. 

\bibitem[Pe]{pe} R. Pellissier, \textit{Cat\'egories enrichies faibles}, Th\`ese,  Universit\'e de Nice-Sophia Antipolis,
June 2002, available at http://math.unice.fr/$^{\sim}$lemaire.

\bibitem[S1]{s1} C. Simpson, \textit{Homotopy over the complex
numbers and generalized cohomology theory}, in \textit{Moduli
of vector bundles (Taniguchi Symposium, December 1994)}, M.
Maruyama ed., Dekker Publ. (1996), 229-263.


\bibitem[S2]{s2} C. Simpson, \textit{The Hodge filtration on non-abelian cohomology}, 
Algebraic geometry---Santa Cruz 1995,  217-281,
Proc. Sympos. Pure Math., 62, Part 2,
Amer. Math. Soc., Providence, RI, 1997. 

\bibitem[S3]{s3} C. Simpson, \textit{Algebraic $n$-stacks}, Preprint  math.AG/9609014.  

\bibitem[S4]{s4} C. Simpson, \textit{A Giraud-type characterisation of the simplicial categories
associated to closed model categories as $\infty$-pretopoi}, Preprint math.AT/9903167.

\bibitem[Tab]{tab} G. Tabuada, \textit{Une structure de cat\'egorie
de mod\`eles de Quillen sur la cat\'egorie des dg-cat\'egories},
Comptes Rendus de l'Acad\'emie de Sciences de Paris, $\mathbf{340}$ (2005), 15-19.

\bibitem[Th]{th} R. Thomas, \textit{A holomorphic Casson invariant for Calabi-Yau 3-folds, 
and bundles on $K3$ fibrations}, 
 J. Differential Geom. $\mathbf{54}$ (2000), no. 2, 367-438.

\bibitem[To1]{to1} B. To\"en, \textit{Homotopical and higher categorical structures in 
algebraic geometry}, habilitation thesis, preprint math.AG/0312262.

\bibitem[To2]{to2} B. To\"en, \textit{The homotopy theory of dg-categories and derived Morita theory}, 
Preprint math.AG/0408337.

\bibitem[To3]{to3} B. To\"en, \textit{Grothendieck rings of higher Artin stacks}, 
Preprint math.AG/0509098.

\bibitem[To4]{to4}B. To\"en, \textit{Derived Hall algebras},
Preprint math.QA/0501343.

\bibitem[To5]{to5}B. To\"en, \textit{Champs affines}, to appear 
in Selecta (also available as Preprint Preprint math.AG/0012219).



\bibitem[To-Va]{tova} B. To\"en, M. Vaqui\'e, \textit{Moduli of objects in dg-categories}, 
preprint math.AG/0503269.

\bibitem[To-Va-Ve]{tvv} B. To\"en, M. Vaqui\'e, G. Vezzosi, \textit{Higher derived categories}, 
in preparation.

\bibitem[HAGI]{hagI} B. To\"en, G. Vezzosi, \textit{Homotopical algebraic geometry I: Topos theory},
Adv. Math. $\mathbf{193}$ (2005), no. 2, 257--372.


\bibitem[HAGII]{hagII} B. To\"en, G. Vezzosi,
\textit{Homotopical algebraic geometry II: Geometric stacks and applications},
to appear in Memoires of the AMS (also available as Preprint math.AG/0404373).

\bibitem[To-Ve1]{tv}  B. To\"en, G. Vezzosi, \textit{Segal topoi and stacks over Segal categories},
December 25, 2002, to appear in Proceedings of the Program \textit{``Stacks, 
Intersection theory and Non-abelian Hodge Theory''},
MSRI, Berkeley, January-May 2002 (also available as Preprint math.AG/0212330).

\bibitem[To-Ve2]{tv2} B. To\"en, G. Vezzosi, \textit{A remark on K-theory and S-categories}, 
Topology $\mathbf{43}$, No. 4 (2004), 765-791.

\bibitem[To-Ve3]{tv3} B. To\"en, G. Vezzosi, \textit{From HAG to DAG: derived moduli stacks}, 
in Axiomatic, enriched and motivic homotopy theory, 
173--216, NATO Sci. Ser. II Math. Phys. Chem., 131, Kluwer Acad. Publ., Dordrecht, 2004.


\end{thebibliography}
\end{document}